\documentclass[12pt]{article}

\usepackage{amssymb, xypic, amsmath, amsthm}

\theoremstyle{plain}
\newtheorem{theorem}{Theorem}[section] 
\newtheorem{corollary}[theorem]{Corollary} 
\newtheorem{lemma}[theorem]{Lemma}
\newtheorem{proposition}[theorem]{Proposition}
\theoremstyle{definition}

\newtheorem{remark}[theorem]{Remark}
\newtheorem{remarks}[theorem]{Remarks}

\numberwithin{equation}{section}

\begin{document}

%
%
%
%
%
%
%
%

\newcommand{\N}{\mathbb{N}}
\newcommand{\Z}{\mathbb{Z}}

\newcommand{\TT}{\mathcal{T}}
\newcommand{\XX}{\mathcal{X}}
\newcommand{\LL}{\mathcal{L}}
\newcommand{\BB}{\mathcal{B}}
\newcommand{\DD}{\mathcal{D}}
\newcommand{\Alg}{{\mathcal{A}}}

\newcommand{\Fin}{\mathsf{Fin}}
\newcommand{\Prim}{{\mathsf{Prim}}}
\newcommand{\type}[1]{\mathsf{type}(#1)}
\newcommand{\rad}{\mathsf{rad}}
\renewcommand{\ker}{\mathsf{ker}}
\renewcommand{\dim}{\mathsf{dim}}
\renewcommand{\min}{\mathsf{min}}
\newcommand{\Char}{\mathsf{char}}

\newcommand{\occursin}{\to}
\newcommand{\ol}{\overline}
\newcommand{\p}{\vdash}
\newcommand{\comp}{\models}
\newcommand{\tensor}{\otimes}
\newcommand{\ass}{\approx}

\newcommand{\setofall}[2]{\mbox{$\{\,#1\,|\,#2\,\}$}}
\newcommand{\tiptopow}[3]{\left\{\begin{array}{ll}
                        #1 & \mbox{ if }#2\\[1mm]
                        #3 & \mbox{ otherwise.}
                        \end{array}\right.}        

\newcommand{\spanofall}[2]{\mbox{$\langle\,#1\,|\,#2\,\rangle$}}

\newcommand{\Setofall}[2]{\mbox{$\Big\{\,#1\,\Big|\,#2\,\Big\}$}}
\newcommand{\tiptop}[4]{\left\{\begin{array}{ll}
                        #1 & #2\\[1mm]
                        #3 & #4
                        \end{array}\right.}        
\newcommand{\length}[1]{\ell(#1)}        

\newcommand{\breakitdown}{%
	\begin{center}$*$\quad $*$\quad $*$\end{center}}

%
%
%
%
%
%
%
%

\allowdisplaybreaks

\begin{large}
  \textbf{%
    The module structure of the Solomon-Tits algebra\\[1mm]
    of the symmetric group}
\end{large}

\bigskip

\textbf{Manfred Schocker}\footnote{%
  supported by Deutsche Forschungsgemeinschaft (DFG Scho-799)}

	\bigskip

  Department of Mathematics\\
  University of Wales Swansea\\
  Singleton Park\\
  Swansea SA2 8PP\\
  United Kingdom\\
	m.schocker@swansea.ac.uk

\bigskip


\bigskip

\today

\medskip

\newpage

\begin{abstract}
  \noindent
  Let $(W,S)$ be a finite Coxeter system.  Tits defined an associative
  product on the set $\Sigma$ of simplices of the associated Coxeter
  complex. The corresponding semigroup algebra is the Solomon-Tits
  algebra of~$W$. It contains the Solomon algebra of~$W$ as the
  algebra of invariants with respect to the natural action of $W$ on
  $\Sigma$. For the symmetric group $S_n$, there is a $1$-$1$
  correspondence between $\Sigma$ and the set of all set compositions
  (or ordered set partitions) of $\{1,\ldots,n\}$. The product
  on $\Sigma$ has a simple combinatorial description in terms
  of set compositions.  
  
  \noindent
  We study in detail the representation theory of the Solomon-Tits
  algebra of $S_n$ over an arbitrary field, and show how our results
  relate to the corresponding results on the Solomon algebra of $S_n$.
  This includes the construction of irreducible and principal
  indecomposable modules, a description of the Cartan invariants, 
  of the Ext-quiver, and of the descending Loewy series.
  Our approach builds on a (twisted) Hopf algebra structure on the
  direct sum of all Solomon-Tits algebras.
\end{abstract}

 \bigskip

\newpage

\section{Introduction} \label{1}
The present work is part of the programme to lift the theory of descent
algebras to the enriched setting of twisted descent algebras as recently
introduced in~\cite{patras-schocker03}. From a different point of
view, this is the attempt to pursue the study of symmetric and
quasi-symmetric functions in non-commuting variables initiated
in~\cite{wolf36} and developed a great deal further
in~\cite{brrz05,bz05,rs05}.

The Solomon-Tits algebra $\TT_n$ of the symmetric group $S_n$ occurs
as a homogeneous component of the free twisted descent algebra. We 
give here a complete description of the module structure of $\TT_n$ over
an arbitrary field. This includes: the construction of primitive
idempotents; a decomposition into principal indecomposable modules; a
description of the Cartan matrix, of the Ext-quiver and, in fact, of
the entire descending Loewy series of $\TT_n$.

\breakitdown

The Solomon algebra $\DD_n$ of $S_n$, discovered by Solomon in the
more general context of finite Coxeter groups, is a subalgebra of the
integral group algebra $\Z S_n$ of $S_n$ and was originally designed
as a noncommutative superstructure of the character ring of the
underlying group~\cite{solomon76}. A huge body of research
papers on the subject, accumulated during the past fifteen years,
provides surprising links to many different fields in geometry,
combinatorics, algebra and topology (\cite{bauer01,bayer-dia92,brown00,erdmann-schocker,garsia89,parisI,joellereut01,parisII,malvenuto-reutenauer95,reutenauer93}, to name but a few; see
\cite{patras-schocker03,schocker-fields} for a more exhaustive list of references).

The Solomon-Tits algebra $\TT_n$ of $S_n$ arises from a semigroup
structure on the set of simplices of the Coxeter complex $\Sigma_n$
associated with $S_n$ and was first considered by Tits in an appendix
to Solomon's original paper~\cite{solomon76}.  The simplices in
$\Sigma_n$ are in $1$-$1$ correspondence with set compositions of
$[n]:=\{1,\ldots,n\}$, that is, with $l$-tuples $(P_1,\ldots,P_l)$ of
mutually disjoint non-empty subsets $P_i$ of $[n]$ such that
$P_1\cup\cdots\cup P_l=[n]$. The associative product on the simplices
in $\Sigma_n$, defined by Tits in geometrical terms, corresponds to
the product $\wedge_n$ on the set $\Pi_n$ of all set
compositions of $[n]$, given by
\begin{eqnarray*}
  \lefteqn{%
    (P_1,\ldots,P_l)\wedge_n (Q_1,\ldots,Q_k)}\\[1mm]
  & := &
  (P_1\cap Q_1,P_1\cap Q_2,\ldots,P_1\cap Q_k,
  \;\ldots\;,
  P_l\cap Q_1,P_l\cap Q_2,\ldots,P_l\cap Q_k)^{\#}
\end{eqnarray*}
for all $(P_1,\ldots,P_l),(Q_1,\ldots,Q_k)\in\Pi_n$
(see~\cite[Section~2]{brown04}).  Here the superscript $\#$ indicates
that empty sets are deleted.  It is readily seen that $P\wedge_n P=P$
for all $P\in\Pi_n$ and that $([n])$ is a two-sided identity in
$(\Pi_n,\wedge_n)$.  This amounts to saying that $(\Pi_n,\wedge_n)$ is
an idempotent semigroup with identity and $\TT_n=\Z \Pi_n$ is the
integral semigroup algebra of this semigroup.

The natural action of $S_n$ on subsets of $[n]$ extends to an action
on $\Pi_n$, defined by
\begin{equation} \label{action}
(P_1,\ldots,P_l)^\pi
:=
(P_1\pi,\ldots,P_l\pi),
\end{equation}
for all $(P_1,\ldots,P_l)\in\Pi_n$ and $\pi\in S_n$. This action
respects the product $\wedge_n$ on $\Pi_n$. Thus the fixed
space $\BB_n$ of $S_n$ in $\TT_n$ is a subalgebra of $\TT_n$.
Bidigare showed that $\BB_n$ is naturally isomorphic to $\DD_n$,
thereby clarifying Tits' original line of reasoning~\cite{bidigare97};
see Section~\ref{3}.

We believe that the often challenging algebraic combinatorics related
to the Solomon algebra can be put to order when $\DD_n$ is viewed as
the ring of invariants in $\TT_n$. The notion of \emph{twisted descent
  algebra} was introduced in~\cite{patras-schocker03} as a formal
framework for this programme.

Let $\N$ denote the set of positive integers. In accordance with the
functorial setup of~\cite{patras-schocker03}, all the semigroups
$(\Pi_A,\wedge_A)$ are considered simultaneously, where $\Pi_A$
consists of the set compositions of an arbitrary finite subset $A$ of 
$\N$ (with the product $\wedge_A$ defined accordingly). The direct sum
$\TT=\bigoplus_A \Z\Pi_A$ carries a second, graded product
$\vee$ and a coproduct $\Delta$, in addition to the internal product
$\wedge_A$ on each homogeneous component, and there are fundamental
rules linking all three structures. The vector space $\TT$ with 
this triple algebraic structure is the free twisted descent algebra.

\breakitdown

In all of what follows, $F$ is an arbitrary field. We will study the
module structure of the Solomon-Tits algebra $(F\Pi_n,\wedge_n)$ over
$F$. Note that any bijection $[n]\to A$ gives rise to an isomorphism
of algebras $F\Pi_n\to F\Pi_A$. In the body of this paper, with an eye
to later applications, results will often be stated for arbitrary 
finite subsets $A$ of $\N$, rather than for $A=[n]$ only.

The representation theory of any finite idempotent semigroup is
governed by its support structure.  The support structure of
$(\Pi_n,\wedge_n)$ arises from the refinement relation $\preceq$ on
$\Pi_n$ defined by $Q\preceq P$ if each component of $Q$ is contained
in a component of $P$. As a general
consequence, the irreducible $F\Pi_n$-modules $M_P$ are one-dimensional 
and naturally labelled by (unordered) set partitions
$P=\{P_1,\ldots,P_l\}$ of $[n]$ or, more conveniently for our purposes,
by those set compositions $P$ contained in
$$
\Pi_n^<
=
\setofall{(P_1,\ldots,P_l)\in\Pi_n}{\min\,P_1<\min\,P_2<\cdots<\min\,P_l}.
$$
This is briefly recovered in Section~\ref{2}. The reader is also
referred to the appendices of~\cite{brown04} as an excellent 
reference on the subject.

In Section~\ref{3} we demonstrate the impact on the Solomon algebra of
these basic observations and, as an illustration of the general concept, prove in a few lines results of Solomon (in characteristic zero) and Atkinson and Willigenburg (in positive characteristic) on the 
radical and the irreducible representations of~$\DD_n$.

The free twisted descent algebra $\TT$ with its additional 
product $\vee$ and its coproduct $\Delta$ is then described in
Section~\ref{4}. Employing the coproduct $\Delta$, allows us 
to construct a
``generic'' primitive idempotent $e_A$ in $(F\Pi_A,\wedge_A)$ for each
finite subset $A$ of $\N$.  As a by-product, we obtain a description
of the primitive elements in $(\TT,\Delta)$ which might be interesting
in its own right.

Each of the products $e_Q=e_{Q_1}\vee\cdots\vee e_{Q_k}$,
where $Q=(Q_1,\ldots,Q_k)\in\Pi_n$, is a primitive idempotent in
$F\Pi_n$, and these elements form a linear basis of $F\Pi_n$. Properly
grouping together basis elements, yields a decomposition of
$F\Pi_n$ into indecomposable left ideals $\Lambda_Q$, where
$Q\in\Pi_n^<$. This is shown in Section~\ref{5}.

Our key result is Theorem~6.2, a multiplication rule for the basis
$\{e_Q\,|\,Q\in\Pi_n\,\}$.  As an immediate consequence, there is a
description of the Cartan invariant $C_{PQ}$ of $F\Pi_n$ (which
describes the multiplicity of $M_P$ in a composition series of
$\Lambda_Q$) as a product of factorials. It also follows 
that $C_{PQ}=0$ unless $Q\preceq P$.

Among other results, in Section~\ref{7}, we show that the occurrence 
of $M_P$ in a composition series of $\Lambda_Q$ is restricted to a
single Loewy layer of $\Lambda_Q$ (if it occurs at all, that is, 
if $Q\preceq P$). As a consequence, the Ext-quiver of
$F\Pi_n$ coincides with the Hasse diagram of the support lattice of
$F\Pi_n$.  This is shown in Section~\ref{8}.

The results on the module structure of $F\Pi_n$ are independent of the
underlying field $F$. We demonstrate in Section~\ref{9} how they 
relate to 
the corresponding results of Garsia and Reutenauer~\cite{garsia-reutenauer89} and Blessenohl and
Laue~\cite{blelau96,blelau02} on the Solomon algebra $\DD_n$
(identified with the subalgebra $\BB_n$ of $\TT_n$) over a field
of characteristic zero. Roughly speaking, all the results 
mentioned above
remain true when the Solomon-Tits algebra is replaced by the Solomon
algebra and compositions $q=(q_1,\ldots,q_k)$ of $n$ are considered
instead of set compositions $Q=(Q_1,\ldots,Q_k)$ of $[n]$.  Most
strikingly, the $k$th radical of $\DD_n$ over $F$ is the invariant
space of the $k$th radical of $\TT_n$ over $F$, for all $k\ge 0$.
Furthermore, via symmetrisation, the primitive idempotents in $\TT_n$
become Lie idempotents in~$\DD_n$.

Finally, we show that the structure of $\TT_n$ as a module for $\DD_n$
is intimately linked to the Garsia-Reutenauer characterisation of
$\DD_n$ as the common stabiliser of a family of vector spaces 
associated with the Poincar\'e-Birkhoff-Witt 
Theorem~\cite[Theorem~4.5]{garsia-reutenauer89}.  More precisely, the
indecomposable $\TT_n$-modules $\Lambda_Q$, when viewed as modules for
$\DD_n$, turn out to be isomorphic to these spaces.

\breakitdown

The representation theory of the modular Solomon algebra (over a field
of positive characteristic $p$) is not yet understood. However, the
little we know has recently proved to be of tremendous help in the
study of modular Lie representations of arbitrary
groups~\cite{bryant-schocker,erdmann-schocker}. It would therefore be
desirable to obtain more information on the modular Solomon algebra.  From the ``twisted'' point of view, the difficulties in
the modular case arise from the multiplicities which occur when
elements of $\TT_n$ are symmetrised.

The Solomon algebra $\DD_W$ of an arbitrary finite Coxeter group $W$
is the algebra of $W$-invariants of the Solomon-Tits algebra
associated with the Coxeter complex of $W$. The approach presented
here for type $A$ may apply to other Coxeter groups as well, to a
certain degree. It could thereby help understand the representation
theory of $\DD_W$ for arbitrary $W$, knowledge of which is very
restricted at the moment.

More generally, there is the question to which extent the results on
the Solo\-mon-Tits algebra of $S_n$ can be lifted to an arbitrary
finite idempotent semigroup. We may, for instance, expect the
Ext-quiver to coincide with the Hasse diagram of the support lattice
for a much larger class of such semigroups.

\section{Support structure and Jacobson radical} \label{2}
Let $S$ be a finite idempotent semigroup with identity.
Due to Brown~\cite{brown04}, the representation theory of $S$
is governed by its support structure. We display
basic definitions and results in this section.

Recall that a partially ordered set $(L,\subseteq)$ is a 
\emph{lattice} if any two elements $x$, $y$ in $L$ have a least 
upper bound $x\vee y$ and a greatest lower bound $x\wedge y$ in
$(L,\subseteq)$.  Assigning to each pair of elements $(x,y)$ in 
$L\times L$ their greatest lower bound $x\wedge y$ defines the 
structure of an abelian idempotent semigroup on $L$.

Associated with the semigroup $S$, there is a \emph{support lattice}
$(L,\subseteq)$ together with a surjective \emph{support map} $s:S\to
L$, such that $(xy)s=xs\wedge ys$, and $xs\subseteq ys$ if and only if
$x=xyx$, for all $x,y\in S$. (We use the
opposite order of that considered in~\cite{brown04} since this seems
more natural for the Solomon-Tits algebra.)

Support lattice and support map of $S$ are unique up to isomorphisms
of ordered sets, thanks to the second condition above.  The first
condition says that $s$ is an epimorphism of semigroups.  

Recall that, if $\Alg$ is a finite-dimensional associative 
$F$-algebra, then the \emph{Jacobson radical} $\rad\,\Alg$ of $\Alg$
is the smallest ideal $R$ of $\Alg$ such that $\Alg/R$ is semisimple.

\begin{theorem}[Brown, 2004] \label{th-brown}
  The support map $s$ extends linearly to an epimorphism of semigroup
  algebras $s:FS\to FL$ such that $\rad\,FS=\ker\,s$.
  
  Furthermore, $FS/\rad\,FS\cong FL$ is abelian and split semisimple,
  that is, isomorphic to a direct sum of copies of the ground field
  $F$.
\end{theorem}

Hence the irreducible representations of $FS$ are all
linear, and in $1$-$1$ correspondence to the elements of the support
lattice $L$.

Combining Brown's theorem with an observation of Bauer~\cite[5.5
Hilfs\-satz]{bauer01} on an arbitrary finite-dimensional algebra with
commutative radical factor, gives the following useful result.

\begin{corollary} \label{cor-brown}
  For any subalgebra $B$ of $FS$, 
  $
  \rad\,B
  =
  B\cap \ker\,s
  $.
  In particular, $B/\rad\,B\cong Bs$ is split semisimple.
\end{corollary}

The short proof follows for the reader's convenience.

\begin{proof}
  Let $R=\rad\,FS$, then $R\cap B$ is a nilpotent ideal of~$B$ and
  therefore contained in $\rad\,B$. Conversely, the only nilpotent
  element of $B/R\cap B\cong (B+R)/R\subseteq FS/R$ is $0$, since
  $FS/R$ is split semisimple by Theorem~\ref{th-brown}.  This implies
  $\rad\,B\subseteq R$ and $ \ker\,s|_B = B\cap\ker\,s = B\cap R =
  \rad\,B $ as asserted.  Furthermore, $FL$ is split semisimple, hence
  so is the subalgebra $Bs$ of $FL$.
\end{proof}

\breakitdown

Let $n\in\N$. The support structure of the semigroup
$(\Pi_n,\wedge_n)$ can be described as follows
(see~\cite[Example~3.2]{brown04}).  Let $Q=(Q_1,\ldots,Q_k)$ and
$P=(P_1,\ldots,P_l)$ be set compositions of $[n]$. We write
$$
Q\preceq P
$$
if each component $Q_i$ of $Q$ is contained in a component $P_j$ of
$P$ or, equivalently, if 
$
Q=Q\wedge_n P\;(=Q\wedge_n P\wedge_n Q)
$.
The relation $\preceq$ on $\Pi_n$ is reflexive and transitive, and
will be of crucial importance in the sequel.  We have $Q\preceq P$ and
$P\preceq Q$ if and only if $Q$ may be obtained by rearranging the
components of $P$. In this case, write
$$
Q\ass P.
$$
The set of equivalence classes, $L$, in $\Pi_n$ with respect to
$\ass$ is a support lattice of $(\Pi_n,\wedge_n)$, with order
inherited from $(\Pi_n,\preceq)$. The greatest lower bound of two
equivalence classes $[P]_\ass$ and $[Q]_\ass$ in $L$ is
$[P]_\ass\wedge [Q]_\ass=[P\wedge_n Q]_\ass$, for all $P,Q\in\Pi_n$.
The support map $s$ sends $P$ to $[P]_\ass$.  Using more illustrative
terms, $L$ may be identified with the set of (unordered) set
partitions of $[n]$ and $s$ with the map that forgets the ordering of
set compositions: $(P_1,\ldots,P_l)s=\{P_1,\ldots,P_l\}$.
Applying Theorem~\ref{th-brown}, we get:

\begin{corollary} \label{abel-set}
  Let $n\in\N$, then $\rad\,F\Pi_n$ is linearly generated by the 
  elements
  $$
  Q-Q'  \quad  (Q,Q'\in\Pi_n,\,Q\ass Q').
  $$
  Furthermore, $F\Pi_n/\rad\,F\Pi_n$ is split semisimple with
  dimension equal to the number of (unordered) set partitions of $[n]$.
\end{corollary}

\section{The Solomon algebra} \label{3}

Let $n\in\N$.  The Solomon algebra $\DD_n$ occurs naturally as a
subalgebra of the Solomon-Tits algebra $\Z\Pi_n$. Thus we can apply
Corollary~\ref{cor-brown} and deduce basic structural information on
$\DD_n$ (over the field $F$).  The results are not new, but their
proofs demonstrate the advantage in viewing the Solomon algebra as the
ring of $S_n$-invariants of the Solomon-Tits algebra at this stage
already. A more detailed analysis of $\DD_n$ follows in Section~\ref{9}.

We recall the necessary definitons.  A finite sequence
$q=(q_1,\ldots,q_k)$ of positive integers with sum $n$ is a
\emph{composition} of $n$, denoted by $q\comp n$.  We write $S_q$ for
the usual embedding of the direct product $S_{q_1}\times\cdots\times
S_{q_k}$ in $S_n$.

The length of a permutation $\pi$ in $S_n$ is the number of inversions
of $\pi$. Each right coset of $S_q$ in $S_n$ contains a unique
permutation of minimal length. Define $\Xi^q$ to be the sum in the
integral group ring $\Z S_n$ of all these minimal coset
representatives of $S_q$ in $S_n$.  Due to Solomon
\cite[Theorem~1]{solomon76}, the $\Z$-linear span of the elements
$\Xi^q$ ($q\comp n$) is a subring of $\Z S_n$ of rank $2^{n-1}$. This
is the \emph{Solomon algebra} $\DD_n$ of~$S_n$.

For any $Q=(Q_1,\ldots,Q_k)\in\Pi_n$, the \emph{type} of $Q$ is
$\type{Q} := (\# Q_1,\ldots,\# Q_k)$.  Hence $\type{Q}\comp n$. Two
set compositions $P,Q\in\Pi_n$ belong to the same $S_n$-orbit if and
only if $\type{P}=\type{Q}$. As a consequence, the ring of
$S_n$-invariants $\BB_n$ in $\Z \Pi_n$ has $\Z$-basis consisting
of the elements
$$
X^q = \sum_{\type{Q}=q} Q \qquad (q\comp n).
$$
Here is Bidigare's remarkable observation (\cite{bidigare97},
see also~\cite[Section~9.6]{brown00}).

\begin{theorem}[Bidigare, 1997] \label{th-bidigare}
  The linear map, defined by $X^q\mapsto \Xi^q$ for all $q\comp n$, is
  an isomorphism of algebras from $\BB_n$ onto $\DD_n$.
\end{theorem}

Now consider the Solomon algebra $\DD_{n,F}=F\tensor_\Z\DD_n$ over
$F$. The results of the previous section yield a description of the
radical of $\DD_{n,F}$ and its quotient in a straightforward way.

Some additional definitions are needed. A composition $r$ of $n$ is a
\emph{partition} of $n$ if $r$ is weakly decreasing. In this case, we
write $r\p n$.  Furthermore, a partition $r$ is \emph{$p$-regular}
(with respect to a positive integer $p$) if no component of $r$ occurs
more than $p-1$ times in $r$.  Finally, we write $q\ass\tilde{q}$ if
$q$ is obtained by rearranging the components of $\tilde{q}$, for all
$q,\tilde{q}\comp n$.

\begin{corollary} \label{rad-Dn}
  The quotient $\DD_{n,F}/\rad\,\DD_{n,F}$ is split semi-simple.
  Furthermore, if $F$ has characteristic zero, then
  $$
  \rad\,\DD_{n,F} 
  =
  \spanofall{\Xi^q-\Xi^{\tilde{q}}}{%
    q,\tilde{q}\comp n,\;q\ass\tilde{q}}_F\,,
  $$
  and the dimension of $\DD_{n,F}/\rad\,\DD_{n,F}$ is equal to the
  number of partitions of $n$.   
  If $F$ has prime characteristic $p$, then
  $$
  \rad\,\DD_{n,F}
  =
  \spanofall{\Xi^q-\Xi^{\tilde{q}}}{%
    q,\tilde{q}\comp n,\;q\ass\tilde{q}}_F 
  \oplus \spanofall{\Xi^r}{r\p  n,\,r\text{ not $p$-regular}}_F\,,
  $$
  and the dimension of $\DD_{n,F}/\rad\,\DD_{n,F}$ is equal to the
  number of $p$-regular partitions of $n$.
\end{corollary}

This is due to Solomon \cite{solomon76} if $F$ has characteristic
zero, and due to Atkinson and Willigenburg \cite{aw97} if $F$ has
prime characteristic.

\begin{proof}
  Theorem~\ref{th-bidigare} allows us to consider
  $\BB_{n,F}=F\tensor_\Z\BB_n$ instead of $\DD_{n,F}$. We will drop
  the index $F$ in what follows.  Let $\tilde{\Pi}_n$ denote the set
  of (unordered) set partitions of $[n]$, and let $s:F\Pi_n\to
  F\tilde{\Pi}_n$ be the linear extension of the support map that
  forgets the ordering.  
  Then by Corollary~\ref{cor-brown}, we have 
  $\rad\,\BB_n=\BB_n\cap\ker\,s$, 
  and $\BB_n/\rad\,\BB_n\cong \BB_ns$ is split semisimple.
  
  For each $q\comp n$, set $c_q:=m_1!\cdots m_n!$, where $m_i$ denotes
  the multiplicity of the entry $i$ in $q$ for all $i\in[n]$.  For
  each $r=(r_1,\ldots,r_l)\p n$, let $x^r$ denote the sum in
  $F\tilde{\Pi}_n$ of all set partitions
  $\tilde{Q}=\{\tilde{Q}_1,\ldots,\tilde{Q}_l\}$ of $[n]$ whose
  elements have cardinalities $r_1,\ldots,r_l$.  Then, if $q\comp n$
  and $r\p n$ with $q\ass r$,
  $$
  X^qs
  =
  \sum_{\type{Q}=q} Qs
  =
  c_q x^r
  =
  c_r x^r.
  $$
  As a consequence, $\BB_ns$ is linearly generated by the elements
  $x^r$ where $r\p n$ such that $\Char\,F$ does not divide $c_r$.
  Furthermore, $(X^q-X^{\tilde{q}})s=0$ whenever $q\ass \tilde{q}$,
  and even $X^qs=0$ whenever $\Char\,F$ divides $c_q$.
  
  Comparing dimensions, completes the proof.
\end{proof}

\section{The free twisted descent algebra} \label{4}
Let $\Fin$ denote the set of all finite subsets of $\N$.  If
$A\in\Fin$ has order $n$, then the set $\Pi_A$ of all set compositions
of $A$ is an idempotent semigroup with identity $(A)$, with respect to
the product defined in the introduction for $A=[n]$.
This product is denoted by $\wedge_A$.  It is clear that the
semigroups $(\Pi_n,\wedge_n)$ and $(\Pi_A,\wedge_A)$ are isomorphic.

Aiming at structural properties of the Solomon-Tits algebra $F\Pi_n$,
it is advantageous to consider several of the semigroups $\Pi_A$
simultaneously and to employ inductive techniques which arise from the
structure of the free twisted descent algebra.

Let $ \Pi=\bigcup_{A\in\Fin} \Pi_A $.  The products $\wedge_A$ extend
to an internal, or \emph{intersection product} $\wedge$ on the direct
sum
$$
F\Pi=\bigoplus_{A\in\Fin} F\Pi_A\,,
$$
by orthogonality: 
$$
P\wedge Q
:=
\tiptop{P\wedge_A Q}{\text{if }\bigcup P=A=\bigcup Q,}{0}{\text{otherwise,}}
$$
for all $P,Q\in\Pi$.  A second, external, or \emph{concatenation
  product} $\vee$ on $F\Pi$ is defined by
$$
P\vee Q
:=
\tiptop{(P_1,\ldots,P_l,Q_1,\ldots,Q_k)}{%
  \text{if }\bigcup P\cap\bigcup Q=\emptyset,}{0}{\text{otherwise,}}
$$
for all $P=(P_1,\ldots,P_l),Q=(Q_1,\ldots,Q_k)\in\Pi$, and
linearity.  If $A=\emptyset$, then $\Pi_A$ consists of the unique set
composition of $A$: the empty tuple $\emptyset$, which acts as a
two-sided identity in $(F\Pi,\vee)$.

If $A\in\Fin$ and $P\in\Pi_A$, set $P|_X:=(P_1\cap X,\ldots,P_l\cap
X)^{\#}\in\Pi_X$ for all $X\subseteq A$.  We define a coproduct on
$F\Pi$ by
$$
\Delta(P)
=
\sum_{X\subseteq A} P|_X\tensor P|_{A\setminus X}
$$
for all $A\in\Fin$, $P\in\Pi_A$, and linearity. This coproduct is
cocommutative.

The algebra $(F\Pi,\vee)$ is a free associative algebra in the
category of $\Fin$-graded vector spaces, with one generator in each
degree, and a free twisted descent algebra (for details, see
\cite{patras-schocker03}).  Background from the theory of twisted
algebras, however, is not needed here.

We recall three simple, but crucial set-theoretical observations
\cite[Theorem~17, Corollary~18]{patras-schocker03}, linking the two
products and the coproduct on $F\Pi$.  Their proofs are sketched for
the reader's convenience. We denote by $\wedge_2$ (respectively, 
by $\vee_2$) the (componentwise) product on $F\Pi\tensor F\Pi$ 
induced by $\wedge$ (respectively, by $\vee$).

\begin{proposition} \label{bi-vee}
  $(F\Pi,\vee,\Delta)$ is a $\Fin$-graded bialgebra, that is:
  if $A,B\in\Fin$ are disjoint, then
  $F\Pi_A\vee F\Pi_B\subseteq F\Pi_{A\cup B}$ and
  $$
  \Delta(f\vee g)
  =
  \Delta(f)\vee_2 \Delta(g)
  $$
  for all $f\in F\Pi_A$, $g\in F\Pi_B$.
\end{proposition}

For the proof, it suffices to consider $P\in \Pi_A$ and $Q\in \Pi_B$,
by linearity. In this case, set $C:=A\cup B$, then
\begin{eqnarray*}
  \Delta(P\vee Q)
  & = &
  \sum_{X\subseteq C}
  (P\vee Q)|_X\tensor (P\vee Q)|_{C\setminus X}\\[1mm]
  & = &
  \sum_{U\subseteq A}
  \sum_{V\subseteq B}
  (P|_U\vee Q|_V)
  \tensor 
  (P|_{A\setminus U}\vee Q|_{B\setminus V})\\[1mm]
  & = &
  \Delta(P)\vee_2 \Delta(Q).\qed
\end{eqnarray*}

\begin{proposition} \label{bi-wedge}
  $(F\Pi,\wedge,\Delta)$ is a bialgebra.
\end{proposition}

The proof is equally simple.\qed

Let $m_\vee:F\Pi\tensor F\Pi\to F\Pi$ be the linearisation of the
product $\vee$ on $F\Pi$.

\begin{proposition}\label{reziprozi}
  $(f\vee g)\wedge h=m_\vee\Big((f\tensor g)\wedge_2 \Delta(h)\Big)$,
  for all $f,g,h\in F\Pi$.
\end{proposition}

Using Sweedler's notation, this reads
$(f\vee g)\wedge h=\sum (f\wedge h^{(1)})\vee (g\wedge h^{(2)})$.

For the proof, it suffices again to consider $P,Q,R\in\Pi$, by
linearity.  Let $A=\bigcup P$, $B=\bigcup Q$ and $C=\bigcup R$, then
$$
m_\vee\Big((P\tensor Q)\wedge_2 \Delta(R)\Big)\\[1mm]
= 
\sum_{X\subseteq C} (P\wedge R|_X)\vee (Q\wedge R|_{C\setminus X}).
$$
This term does not vanish if and only if $A\subseteq C$ and
$B=C\backslash A$, that is, if $C$ is the disjoint union of $A$ and
$B$.  The same is true for the term $(P\vee Q)\wedge R$.  And in this
case,
\begin{eqnarray*}
  \lefteqn{%
    m_\vee\Big((P\tensor Q)\wedge_2 \Delta(R)\Big)}\\[3mm]
  & = &
  (P\wedge R|_A)\vee (Q\wedge R|_B)\\[3mm]
  & = &
  (P_1\cap R_1,\ldots,P_1\cap R_m,
  \;\ldots\;,
  P_l\cap R_1,\ldots,P_l\cap R_m)^{\#}\\[2mm]
  &&
  \mbox{\hphantom{$(P\wedge R|_A)$}}
  \vee
  (Q_1\cap R_1,\ldots,Q_1\cap R_m,
  \;\ldots\;,
  Q_k\cap R_1,\ldots,Q_k\cap R_m)^{\#}\\[3mm]
  & = &
  (P\vee Q)\wedge R,
\end{eqnarray*}
where $P=(P_1,\ldots,P_l)$, $Q=(Q_1,\ldots,Q_k)$ and
$R=(R_1,\ldots,R_m)$.\qed

Proposition~\ref{reziprozi} often allows us to transfer calculations 
from $(F\Pi_A,\wedge_A)$ to $(F\Pi,\vee,\Delta)$ and $(F\Pi_B,\wedge_B)$
for some (proper) subsets $B$ of $A$. This inductive idea will be
frequently used in what follows.

\breakitdown

Let $A\in\Fin$. An element $e$ of $F\Pi_A$ is an \emph{idempotent} if
$e^2=e\wedge e=e$ and $e\neq 0$. Such an idempotent is
\emph{primitive} if the left ideal $F\Pi_A\wedge e$ of $F\Pi_A$ is
indecomposable as an $F\Pi_A$-module.  (Equivalently, whenever $f,g\in
F\Pi_A$ such that $e=f+g$, $f^2=f$, $g^2=g$ and $f\wedge
g=0=g\wedge f$, then $f=0$ or $g=0$.)

An element $E\in F\Pi$ is \emph{$\Delta$-primitive} if
$\Delta(E)=E\tensor\emptyset+\emptyset\tensor E$.

In concluding this section, we illustrate the inductive method by
exploring a relation between the $\Delta$-primitive elements of
$F\Pi_A$ and certain primitive idempotents in $F\Pi_A$.

\begin{corollary} \label{cor-rp}
  Let $A\in\Fin$ and suppose that $E=\sum_{P\in\Pi_A} k_PP\in F\Pi_A$
  is $\Delta$-primitive, then
  $
  P\wedge E=0
  $
  for all $P\in\Pi_A\backslash\{(A)\}$. 

  In particular, $E\wedge E=k_{(A)}E$.
\end{corollary}

\begin{proof}
  If $P\in\Pi_A\backslash\{(A)\}$, say, $P=R\vee S$ with
  $R,S\in\Pi\setminus\{\emptyset\}$, then
  $$
  P\wedge E
  =
  (R\vee S)\wedge E
  =
  m_\vee\Big(
  (R\tensor S)\wedge_2 (E\tensor\emptyset+\emptyset\tensor E)
  \Big)
  =
  0,
  $$
  by Proposition~\ref{reziprozi}.
  It follows that
  $E\wedge E=k_{(A)}(A)\wedge E=k_{(A)}E$. 
\end{proof}

Note that, if $k_{(A)}\neq 0$, then $e=1/k_{(A)} E$ is in fact a
\emph{primitive} idempotent in $F\Pi_A$.  For, if $f,g\in F\Pi_A$ such
that $e=f+g$, $f^2=f$, $g^2=g$ and $f\wedge g=0=g\wedge f$, then
$e\wedge f=f$ and $e\wedge g=g$. Hence $f$ and $g$ are
$\Delta$-primitive as well, by Proposition~\ref{bi-wedge}.  The
preceding corollary implies $c_{(A)}g=f\wedge g=0$ and
$c_{(A)}f=f\wedge f=f$, where $c_{(A)}$ denotes the coefficient of
$(A)$ in $f$. Thus $f=0$ or $g=0$.

If $P=(P_1,\ldots,P_l)\in\Pi$, then $\length{P}:=l$ is the
\emph{length} of $P$.  We set
$$
\Pi_A^*
:=
\setofall{(P_1,\ldots,P_l)\in\Pi_A}{\min\, A\in P_1}
$$
for all $A\in\Fin$.

\begin{lemma}\label{primitiv}
  Let $A\in\Fin$, then the element
  $$
  e_A = \sum_{P\in\Pi_A^*} (-1)^{\length{P}-1} P
  $$
  is $\Delta$-primitive.  In particular, $e_A$ is a primitive
  idempotent in $F\Pi_A$.
\end{lemma}

\begin{proof}
  Let $R,S\in\Pi\setminus\{\emptyset\}$ such that $R\vee S\in\Pi_A$.
  We need to show that the coefficient $ c_{R,S} $ of $R\tensor S$ in
  $\Delta(e_A)$ is zero.  Let $X=\bigcup R$. It suffices to consider
  the case where $a^*:=\min\,A\in X$, since $\Delta$ is cocommutative.
  We have
  $$
  c_{R,S} = \sum_{P\in \XX}(-1)^{\length{P}-1},
  $$
  where $\XX=\setofall{P\in\Pi_A^*}{P|_X=R,\,P|_{A\setminus X}=S}$.
  Let $P=(P_1,\ldots,P_l)\in\XX$, then there exists an index $i\in[l]$
  such that $P_i\setminus X\neq\emptyset$, since $X\neq A$.  Choose
  $i$ minimal with this property.  If also $P_i\cap X\neq\emptyset$,
  then set
  $$
  P' 
  := 
  (P_1,\ldots,P_{i-1},P_i\cap X,P_i\setminus X,P_{i+1},\ldots,P_l)
  \in 
  \XX.
  $$
  If $P_i\cap X=\emptyset$, then $a^*\in P_1\cap X$ implies that 
  $i>1$, and we define
  $$
  P'
  :=
  (P_1,\ldots,P_{i-2},P_{i-1}\cup P_i,P_{i+1},\ldots,P_l)
  \in
  \XX.
  $$
  Then $(P')'=P$ and $(-1)^{\length{P'}}=-(-1)^{\length{P}}$ for 
  all $P\in\XX$,
  that is, $P\mapsto P'$ defines is a sign-reversing pairing of the 
  summands of $c_{R,S}$.
  This implies that $c_{R,S}=0$, hence that $e_A$ is
  $\Delta$-primitive.  The additional claim follows from
  Corollary~\ref{cor-rp} and its subsequent remark.
\end{proof}

The idempotent $e_A$ is displayed in Table~1 for some
particular choices of $A$. Note that several curly brackets and commas
have been omitted; for instance, $(12,3)$ means $(\{1,2\},\{3\})$.
\begin{table}[htbp]
	$$
	\begin{array}{lcl}
	  e_{\{1\}}   & = & (1),\\[1mm]
	 	e_{\{1,2\}}  & = & (12) - (1,2),\\[1mm]
		e_{\{1,2,3\}} 
		& = & 
		(123)-(12,3)-(1,23)-(13,2)+(1,2,3)+(1,3,2),\\[1mm]
		e_{\{2,5,6\}} 
		& = & 
		(256)-(25,6)-(2,56)-(26,5)+(2,5,6)+(2,6,5).
	\end{array}
	$$
  \caption{The idempotent $e_A$ for $A=\{1\}$, $\{1,2\}$, $\{1,2,3\}$ 
  			and $\{2,5,6\}$.}
  \label{tabx}
\end{table}

The Lie product $\circ$ associated with the product $\vee$ on $F\Pi$
is defined by
$$
f\circ g
=
f\vee g - g\vee f
$$
for all $f,g\in F\Pi$. The set $\Prim\,F\Pi$ of all 
$\Delta$-primitive elements in $F\Pi$ is a Lie subalgebra of
$(F\Pi,\circ)$.  As a by-product of the above considerations, 
there is the following description of $\Prim\,F\Pi$.

\begin{corollary} \label{primlie}
  The $\Delta$-primitive Lie algebra of $F\Pi$ is
  $$
  \Prim\,F\Pi=\bigoplus_{A\in\Fin} e_A\wedge F\Pi_A\,.
  $$
  Its $A$-graded component $e_A\wedge F\Pi_A$ is equal to the
  linear span of all idempotents $e\in F\Pi_A$ such that
  $e\wedge F\Pi_A=e_A\wedge F\Pi_A$, and has dimension $2\#\Pi_{n-1}$.
\end{corollary}

\begin{proof}
  Let $A\in\Fin$, and suppose $e\in F\Pi_A$ is $\Delta$-primitive.
  Then $e=(A)\wedge e=e_A\wedge e\in e_A\wedge F\Pi_A$, by
  Corollary~\ref{cor-rp}.  Conversely, $e_A\wedge f$ is
  $\Delta$-primitive for any $f\in F\Pi_A$, by
  Proposition~\ref{bi-wedge}. 
  
  The set of idempotents $e$ in $F\Pi_A$ with $e\wedge
  F\Pi_A=e_A\wedge F\Pi_A$ is 
  $$
  e_A+e_A\wedge F\Pi_A\wedge ((A)-e_A),
  $$
  as a general fact.  The linear span of this set is equal to
  $e_A\wedge F\Pi_A$, as claimed, since $P\wedge e_A=0$ for all
  $P\in\Pi_A$ with $\length{P}>1$.
  
  The dimension formula will be obtained in
  Remark~\ref{remarks-cartan}~(1).
\end{proof}

\section{Primitive idempotents and principal indecomposable modules}
\label{5}
We are now in a position to study the module structure of the
Solomon-Tits algebra in more detail. In this section, for any finite
subset $A$ of $\N$, we construct a linear basis of $F\Pi_A$ which
consists of primitive idempotents. Properly grouping together basis
elements, we then obtain a decomposition of $F\Pi_A$ into
indecomposable left ideals.

In what follows, $E_A$ is a $\Delta$-primitive idempotent 
in $F\Pi_A$, that is, we have
$\Delta(E_A)=E_A\tensor\emptyset+\emptyset\tensor E_A$
and 
$E_A\wedge E_A=E_A$,
for all $A\in\Fin$. (One possible choice for $E_A$ is the element
$e_A$, defined in Lem\-ma~\ref{primitiv}. We shall make other choices
for $E_A$ in Section~\ref{9}.) Note that
$E_A\in (A)+\langle\,\Pi_A\backslash\{(A)\}\,\rangle_F$ by Corollary~\ref{cor-rp}.
We set
$$ 
E_Q := E_{Q_1}\vee\cdots\vee E_{Q_k}
$$
for all $Q=(Q_1,\ldots,Q_k)\in\Pi$. Then
\begin{equation}
  \label{triangular}
  E_Q\in Q+\spanofall{R\in\Pi_A}{R\preceq Q,\,R\not\ass Q}_F
\end{equation}
for all $A\in\Fin$, $Q\in\Pi_A$. Here $R\preceq Q$ means that each
component of $R$ is contained in a component of $Q$, and $R\ass Q$ 
means that $R$ is a rearrangement of $Q$ (as in the case $A=[n]$).  
From~\eqref{triangular}, using a triangularity argument, we get:

\begin{proposition} \label{E-basis}
  $\setofall{E_Q}{Q\in\Pi_A}$ is a linear basis of $F\Pi_A$, for all
  $A\in\Fin$.
\end{proposition}

The basis $\setofall{e_Q}{Q\in\Pi_3}$ of $F\Pi_3$ (arising
from the idempotents $e_A$ defined in Lemma~\ref{primitiv}), is
displayed in Table~2.
\begin{table}[htbp]
	$$
  \begin{array}{lll@{\hspace*{10ex}}lll}
    e_{(1,2,3)} & = & (1,2,3),
    &
    e_{(1,3,2)} & = & (1,3,2),\\[1mm]
    e_{(2,1,3)} & = & (2,1,3),
    &
    e_{(2,3,1)} & = & (2,3,1),\\[1mm]
    e_{(3,1,2)} & = & (3,1,2),
    &
    e_{(3,2,1)} & = & (3,2,1),\\[1mm]
    e_{(12,3)} & = & (12,3)-(1,2,3),
    &
    e_{(3,12)} & = & (3,12)-(3,1,2),\\[1mm]
    e_{(1,23)} & = & (1,23)-(1,2,3),
    &
    e_{(23,1)} & = & (23,1)-(2,3,1),\\[1mm]
    e_{(13,2)} & = & (13,2)-(1,3,2),
    &
    e_{(2,13)} & = & (2,13)-(2,1,3),\\[1mm]
    e_{(123)}  & = & 
    \multicolumn{4}{l}{(123)-(12,3)-(1,23)-(13,2)+(1,2,3)+(1,3,2).}
  \end{array}
  $$
  \caption{A basis of $F\Pi_3$ consisting of primitive idempotents.}
  \label{tab0}
\end{table}
Some technical preparations are needed to describe the left-regular
representation of $F\Pi_A$ in terms of this basis. Let $Q=(Q_1,\ldots,Q_k)\in \Pi_A$. We set
$$
Q_I:=(Q_{i_1},\ldots,Q_{i_m})
$$
for any subset $I=\{i_1,\ldots,i_m\}$ of $[k]$ with
$i_1<\cdots<i_m$.  If $P=(P_1,\ldots,P_l)\in \Pi_A$ such that
$Q\preceq P$, then there exists a unique set composition
$I=(I_1,\ldots,I_l)$ of the index set $[k]$ of $Q$ such that
$Q_{I_j}\in \Pi_{P_j}$ for all $j\in[l]$. Note that, in this case,
$P\wedge Q=Q_{I_1}\vee\cdots\vee Q_{I_l}$. For example, if
$A=\{3,5,6,7,8\}$, $Q=(8,67,3,5)$ and $P=(367,58)$, then $Q\preceq P$,
$I=(23,14)$ and $P\wedge Q=(67,3,8,5)$.

\begin{lemma} \label{mult-red}
  Let $A\in\Fin$ and $P=(P_1,\ldots,P_l),Q=(Q_1,\ldots,Q_k)\in\Pi_A$.
  Suppose $f^{(j)}\in F\Pi_{P_j}$ for all $j\in[l]$, then
  $$
  (f^{(1)}\vee\cdots\vee f^{(l)})\wedge E_Q
  =
  \tiptop{%
    (f^{(1)}\wedge E_{Q_{I_1}})\vee\cdots\vee(f^{(l)}\wedge E_{Q_{I_l}})}
  {\text{if }Q\preceq P,}{0}{\text{otherwise,}}
  $$
  where, in case $Q\preceq P$, $(I_1,\ldots,I_l)$ is the unique set
  composition in $\Pi_k$ such that $Q_{I_j}\in\Pi_{P_j}$ for all
  $j\in[l]$.
\end{lemma}

\begin{proof}
  For $l=1$, there is nothing to prove (since $Q\preceq (A)$). 
  Let $l>1$, then Proposition~\ref{bi-vee} implies that
  $$
  \Delta(E_Q)
  =
  \Delta(E_{Q_1})\vee_2\cdots\vee_2\Delta(E_{Q_k})
  =
  \sum_{I\subseteq[k]} E_{Q_I}\tensor E_{Q_{[k]\setminus I}}\,.
  $$
  Thus, setting $\tilde{f}:=f^{(2)}\vee\cdots\vee f^{(l)}$
  and applying Proposition~\ref{reziprozi}, we get
  $$
  (f^{(1)}\vee\tilde{f})\wedge E_Q
  =
  \sum_{I\subseteq[k]} 
  (f^{(1)}\wedge E_{Q_I})\vee(\tilde{f}\wedge E_{Q_{[k]\setminus I}}).
  $$
  This term vanishes unless there exists a subset $I_1$ of $[k]$
  such that $Q_{I_1}\in\Pi_{P_1}$, in which case $
  (f^{(1)}\vee\tilde{f})\wedge E_Q = (f^{(1)}\wedge
  E_{Q_{I_1}})\vee(\tilde{f}\wedge E_{Q_{[k]\setminus I_1}}) $. The
  proof may be completed by a simple inductive step.
\end{proof}

Special cases of the preceding lemma are:
\begin{equation}
  \label{left-mod}
  P\wedge E_Q
  =
  \tiptop{E_{P\wedge Q}}{\text{if }Q\preceq P,}{0}{\text{otherwise,}}
  \mbox{ for all $P,Q\in\Pi_A$}
\end{equation}
(with $f^{(i)}=(P_i)$ for all $i\in[l]$), and 
\begin{equation}
  \label{prim-idem}
  E_P\wedge E_Q
  =
  E_P
  \mbox{ whenever $P\ass Q$}
\end{equation}
(with $f^{(i)}=E_{P_i}$ for all $i\in[l]$).
Considering $P=Q$, we obtain:

\begin{corollary} \label{idem-cor}
  $E_Q\wedge E_Q=E_Q$ for all $Q\in \Pi$.
\end{corollary}

As we shall see below, each $E_Q$ is in fact a \emph{primitive}
idempotent.
Note that \eqref{triangular} and \eqref{left-mod} imply furthermore
that
\begin{equation}
  \label{mult-triangle}
  E_P\wedge E_Q=0\mbox{ unless }Q\preceq P.
\end{equation}

\breakitdown

In what follows, all modules are left modules.  If $\Alg$ is a
finite-dimensional associative algebra with identity and $M$ is 
an $\Alg$-module, then the \emph{$\Alg$-radical} $\rad_\Alg M$ 
of $M$ is the intersection of all
maximal $\Alg$-submodules of $M$.  In particular,
$\rad\,\Alg=\rad_\Alg\Alg$ is the Jacobson radical of $\Alg$, and
$\rad_\Alg M=(\rad\,\Alg)\,M$.

We observe that the set
$$
\Pi_A^<
:=
\setofall{(Q_1,\ldots,Q_k)\in\Pi_A}{\min\,Q_1<\cdots<\min\,Q_k}
$$
is a transversal for the equivalence classes in $\Pi_A$ with
respect to $\ass$, for each $A\in\Fin$. (Hence the irreducible
$F\Pi_A$-modules are in $1$-$1$ correspondence to the elements of
$\Pi_A^<$, by Corollary~\ref{abel-set}.)

\begin{theorem} \label{indecs}
  For each $Q\in\Pi_A$,
  $$
  \Lambda_Q
  :=
  F\Pi_A\wedge E_Q
  =
  \spanofall{E_{Q'}}{Q'\in\Pi_A,\;Q'\ass Q}_F
  $$
  is an indecomposable $F\Pi_A$-left module with
  radical
  $$
  \rad_{F\Pi_A}\Lambda_Q
  =
  \spanofall{E_Q-E_{Q'}}{Q'\in\Pi_A,\;Q'\ass Q}_F
  $$
  of codimension $1$.  In particular,
  $$
  F\Pi_A
  =
  \bigoplus_{Q\in\Pi_A^<} \Lambda_Q
  $$
  is a decomposition into indecomposable submodules, and
  $$
  \rad\,F\Pi_A
  =
  \spanofall{E_Q-E_{Q'}}{Q\in\Pi_A^<,\;Q'\in\Pi_A,\;Q'\ass Q}_F\,.
  $$
\end{theorem}

\begin{proof}
  The elements $E_{Q'}$, $Q'\ass Q$, constitute a linear basis of
  $\Lambda_Q=F\Pi_A\wedge E_Q$, by \eqref{left-mod}, \eqref{prim-idem}
  and Proposition~\ref{E-basis}.  In particular,
  $\Lambda_Q=\Lambda_{Q'}$ whenever $Q\ass Q'$. Hence $F\Pi_A =
  \bigoplus_{Q\in\Pi_A^<} \Lambda_Q$ is a decomposition into left
  ideals by Proposition~\ref{E-basis}.
  
  Since the dimension of $F\Pi_A/\rad\,F\Pi_A$ is equal to $\#\Pi_A^<$,
  by Corollary~\ref{abel-set}, it follows that
  $\Lambda_Q/\rad_{F\Pi_A}\Lambda_Q$ is one-dimensional for all $Q$.
  This implies that $\Lambda_Q$ is indecomposable and also the
  description of $\rad_{F\Pi_A}\Lambda_Q$.
\end{proof}

As an immediate consequence, each of the idempotents $E_Q$ is
primitive in $F\Pi_A$. Besides,
$\dim\,\Lambda_Q=\length{Q}!$ for all $Q\in\Pi_A$.

\begin{corollary}
  The one-dimensional spaces 
  $
  M_Q=\Lambda_Q/\rad_{F\Pi_A}\Lambda_Q\,
  $,
  $Q\in\Pi_A^<$, form a complete set of mutually non-isomorphic
  irreducible $F\Pi_A$-modules.
\end{corollary}

\begin{remark} \label{idem-choice} 
	Suppose $\tilde{E}_A$ is another $\Delta$-primitive idempotent
	in $F\Pi_A$, for all $A\in\Fin$, and 
	$\tilde{\Lambda}_Q=F\Pi_A\wedge\tilde{E}_Q$ 
	denotes the corresponding indecomposable $F\Pi_A$-module
	for all $Q\in\Pi_A$. Then $\Lambda_Q\cong \tilde{\Lambda}_Q$
	as $F\Pi_A$-modules.
	
	Indeed, from Corollary~\ref{cor-rp} we get
	that $E_A\wedge\tilde{E}_A=\tilde{E}_A$ for all $A\in\Fin$.
	Hence $E_{Q'}\wedge \tilde{E}_Q=\tilde{E}_{Q'}$
	whenever $Q'\ass Q$ by Lemma~\ref{mult-red}. It follows
	that $f\mapsto f\wedge \tilde{E}_Q$ defines an isomorphism
	from $\Lambda_Q$ onto $\Lambda_{\tilde{Q}}$. 
\end{remark}

\section{Cartan invariants} \label{6}
Taking $E_A=e_A$ for all $A\in\Fin$, we obtain the linear basis
$\setofall{e_Q}{Q\in\Pi_A}$ of $F\Pi_A$ from
Lemma~\ref{primitiv} and Proposition~\ref{E-basis}.
It is well adapted to the module structure of
$F\Pi_A$, as will be further stressed in the sections that follow.

The key result is Theorem~\ref{multiformel}, a multiplication rule
for the basis $\{e_Q\,|\,Q\in\Pi_A\}$, which allows us to
determine the Cartan matrix of $F\Pi_A$ at once.

The Lie product $\circ$ associated with the product $\vee$ on $F\Pi$
has been considered at the end of Section~\ref{4} already. It occurs in a
natural way when two basis elements $e_P$ and $e_Q$ are multiplied
together.  For example, if $P=(A)$ and $Q=(X,Y)\in\Pi_A$, then
\begin{equation} \label{multi-l2}
  e_A\wedge e_{(X,Y)}
  =
  \tiptopow{e_X\circ e_Y}{\min\,A\in Y,}{0}
\end{equation}
Indeed, if $\min\,A\in X$, then
$$
e_A\wedge e_{(X,Y)}
=
\sum_{P\in\Pi_A^*} (-1)^{\length{P}-1} P\wedge e_{(X,Y)}
=
(A)\wedge e_{(X,Y)}-(X,Y)\wedge e_{(X,Y)}
=
0,
$$
by \eqref{left-mod}. Similarly, if $\min\,A\in Y$, then
$$
e_A\wedge e_{(X,Y)}
=
(A)\wedge e_{(X,Y)}-(Y,X)\wedge e_{(X,Y)}
=
e_{(X,Y)}-e_{(Y,X)}.
$$
A more systematical approach follows.  The (right-normed)
\emph{Dynkin mapping} $F\Pi\to F\Pi,\,f\mapsto f^\circ$ is defined
recursively by $ (A)^\circ=(A) $ for all $A\in\Fin$,
$$
(Q_1,\ldots,Q_k)^\circ
:=
Q_1\circ\Big((Q_2,\ldots,Q_k)^\circ\Big),
$$
for all $(Q_1,\ldots,Q_k)\in\Pi$ with $k>1$, and linearity.  
For $A\in\Fin$, set
$$
\Pi_A^\dagger:=\setofall{(Q_1,\ldots,Q_k)\in\Pi_A}{\min\,A\in Q_k}.
$$
We will need the following folklore result on right-normed
multilinear Lie monomials (see \cite[Proposition~2.3]{schocker-rel}
for the left-normed case).

\begin{proposition} \label{lie-help}
  Let $A\in\Fin$ and $Q=(Q_1,\ldots,Q_k)\in\Pi_A^\dagger$, then
  \begin{equation} \label{r-n-tri}
  Q^\circ
  \in
  Q+\spanofall{R\in\Pi_A}{R\ass Q,\,R\notin \Pi_A^\dagger}_F\,.
  \end{equation}
  In particular,  
  $\setofall{\tilde{Q}^\circ}{\tilde{Q}\in\Pi_A^\dagger,\,Q\ass \tilde{Q}}$
  is a linear basis of
  $\spanofall{R^\circ}{R\ass Q}_F$ of order $(k-1)!$.
\end{proposition}

\begin{proof}
	The first claim follows by a simple induction on $k$ and implies
	linear independency of the set considered in the second claim.  The
	rest follows by comparing dimensions.
\end{proof}

For all $Q=(Q_1,\ldots,Q_k),P=(P_1,\ldots,P_l)\in\Pi_A$, we write
$$
Q\preceq^\dagger P
$$
if there exists a set composition $(I_1,\ldots,I_l)\in\Pi_{[k]}$
such that $Q_{I_j}\in\Pi_{P_j}^\dagger$ for all $j\in[l]$.  In this
case, $(I_1,\ldots,I_l)$ is unique, and certainly $Q\preceq P$.  For
example, $(4,5,13,2)\preceq(123,45)$, but \emph{not}
$(4,5,13,2)\preceq^\dagger(123,45)$, since $I_1=\{3,4\}$
and $Q_{I_1}=(13,2)\notin\Pi_{\{1,2,3\}}^\dagger$.

\begin{theorem} \label{multiformel}
  Let $P=(P_1,\ldots,P_l),Q=(Q_1,\ldots,Q_k)\in\Pi_A$, then
  $$
  e_P\wedge e_Q
  =
  \tiptop{%
    e_{Q_{I_1}^\circ}\vee\cdots\vee e_{Q_{I_l}^\circ}}{%
    \text{if $Q\preceq^\dagger P$,}}{%
    0}{\text{otherwise}}
  $$
  where, in case $Q\preceq^\dagger P$, the set composition
  $(I_1,\ldots,I_l)$ of $[k]$ is so chosen that
  $Q_{I_j}\in\Pi_{P_j}^\dagger$ for all $j\in[l]$.
  
  In particular, $e_A\wedge e_Q=e_{Q^\circ}$ if $Q\in\Pi_A^\dagger$,
  and $e_A\wedge e_Q=0$ otherwise.
\end{theorem}

Here the map $e:Q\mapsto e_Q$ has been extended linearly to
$F\Pi$, so that by definition
$$
e_{\raisebox{-4pt}{$\Sigma f_RR$}}
=
\sum f_R\,e_R
$$
for all $\sum f_RR\in F\Pi$.  For example, we have
$e_{(123,45)}\wedge e_{(4,5,13,2)}=0$, since
$(4,5,13,2)\preceq^\dagger (123,45)$ does not hold, while
$$
e_{(234,1)}\wedge e_{(34,1,2)}
=
e_{(34,2)^\circ}\vee e_{(1)}
=
e_{(34,2)-(2,34)}\vee e_{(1)}
=
e_{(34,2,1)}-e_{(2,34,1)}\,.
$$

\begin{proof}[Proof of Theorem~\ref{multiformel}.]
  The goal is to prove the assertion in the special case where $P=(A)$,
  mentioned in the second part, since the general case then follows
  from Lemma~\ref{mult-red}, applied to $f^{(i)}=e_{P_i}$ 
  for all $i\in[l]$.
  
  This will be done in four steps. Choose $i\in[k]$ such that
  $a^*:=\min\,A\in Q_i$.
  
  \emph{Step 1}.  Suppose $P\in\Pi_A^*$ such that $Q\preceq P$ and set
  $Q':=P\wedge Q$. Then $a^*\in P_1$ and $Q'\ass Q$. Furthermore, from
  $Q'=(Q'_1,\ldots,Q'_k)=(P_1\cap Q_1,\ldots,P_1\cap Q_i,\ldots)^{\#}$
  it follows that $a^*\in Q'_j$ for some $j\le i$; and we have $i=j$
  if and only if $\bigcup_{m\le i} Q_m\subseteq P_1$.
  
  \emph{Step 2}.  Suppose that $i<k$, and let $Q'\ass Q$ such that
  $a^*\in Q'_i$. Put 
  $\XX:=\setofall{P\in\Pi_A^*}{Q'=P\wedge Q}$.
  Then $\bigcup_{m\le i} Q_m\subseteq P_1$ for each $P\in\XX$, by
  Step~1.  Hence the set $\XX$ decomposes into subsets $\XX_=$ 
  and $\XX_{\neq}$, where $P=(P_1,\ldots,P_l)\in\XX$ belongs to 
  $\XX_=$ or $\XX_{\neq}$ according as $\bigcup_{m\le i} Q_m= P_1$ 
  or not. 
  Note that $i<k$ implies that $\length{P}>1$ for all $P\in\XX_=$.
  The mapping
  $$
  \XX_=\to\XX_{\neq},\,
  P\mapsto\tilde{P}
  :=
  (P_1\cup P_2,P_3,\ldots,P_l)
  \mbox{\hphantom{$\displaystyle%
      \bigcup_{m\le i} Q_m,\tilde{P}_1\setminus\bigcup_{m\le i} Q_m$}}
  $$
  is a bijection with inverse given by
  $$
  \XX_{\neq}\to\XX_=,\,
  \tilde{P}\mapsto P
  :=
  (\bigcup_{m\le i} Q_m,
   \tilde{P}_1\setminus\bigcup_{m\le i} Q_m,
   \tilde{P}_2,\ldots,\tilde{P}_l),
  \mbox{\hphantom{$P_1\cup P_2$}}
  $$
  and $\length{\tilde{P}}=\length{P}-1$ for all $P\in \XX_=$.
  
  \emph{Step 3}. We are now in a position to prove $e_A\wedge e_Q=0$ if
  $i<k$, by induction on $i$. For, 
  \begin{eqnarray*}
    e_A\wedge e_Q
    & = &
    e_A\wedge(e_A\wedge e_Q),\mbox{ by Lemma~\ref{primitiv}}\\[3mm]
    & = &
    e_A\wedge
    \sum_{Q'\ass Q}
     \sum_{\substack{P\in\Pi_A^*\\[2pt] Q'=P\wedge Q}} 
      (-1)^{\length{P}-1} e_{Q'},
      \mbox{ by~\eqref{left-mod}}\\[1mm]
    & = &
    e_A\wedge
    \sum_{j=1}^i
    \sum_{\substack{Q'\ass Q\\[2pt] a^*\in Q'_j}}
    \Big(
    \sum_{\substack{P\in\Pi_A^*\\[2pt] Q'=P\wedge Q}} (-1)^{\length{P}-1}\Big) 
    e_{Q'},
      \mbox{ by Step~1}\\[1mm]
    & = &
    \sum_{j=1}^{i-1}
     \sum_{\substack{Q'\ass Q\\[2pt] a^*\in Q'_j}}
    \Big(
      \sum_{\substack{P\in\Pi_A^*\\[2pt] Q'=P\wedge Q}}
       (-1)^{\length{P}-1}
    	\Big) 
    e_A\wedge e_{Q'}\,,\mbox{ by Step~2.}
  \end{eqnarray*}
  Thus, if $i=1$, then $e_A\wedge e_Q=0$ follows directly, while for
  $i>1$, we may conclude by induction.
  
  \emph{Step 4}. Assume now that $i=k$. We will show that 
  $e_{Q^\circ}=e_A\wedge e_Q$ by induction on $k=\length{Q}$.
  
  For $k=1$, this is Lemma~\ref{primitiv}.  Let $k>1$, and set
  $X=Q_1$, $\tilde{Q}=(Q_2,\ldots,Q_k)$ and $Y=\bigcup \tilde{Q}$,
  then
  \begin{eqnarray*}
    e_{Q^\circ}
    & = &
    e_X\vee e_{\tilde{Q}^\circ}-e_{\tilde{Q}^\circ}\vee e_X\\[3mm]
    & = &
    e_X\vee (e_Y\wedge e_{\tilde{Q}})
    -
    (e_Y\wedge e_{\tilde{Q}})\vee e_X,
    \mbox{ by induction}\\[3mm]
    & = &
    (e_{(X,Y)}-e_{(Y,X)})\wedge e_Q,\mbox{ by Lemma~\ref{mult-red}}\\[3mm]
    & = &
    e_A\wedge e_{(X,Y)}\wedge e_Q,
    \mbox{ by~\eqref{multi-l2}}\\[3mm]
    & = &
    e_A\wedge \Big(e_X\vee (e_Y\wedge e_{\tilde{Q}})\Big),
    \mbox{ by Lemma~\ref{mult-red}}\\[3mm]
    & = &
    e_A\wedge(e_X\vee e_{\tilde{Q}^\circ}),
    \mbox{ by induction}\\[3mm]
    & = &
    e_A\wedge e_Q,
    \mbox{ by~\eqref{r-n-tri} and Step~3.}
  \end{eqnarray*}
  The proof is complete.
\end{proof}

\begin{corollary} \label{complement}
  The elements $e_T$, $T\in \Pi_A^<$, form a complete set of mutually
  orthogonal primitive idempotents in $F\Pi_A$, and
  $\sum_{T\in\Pi_A^<} e_T=(A)$. 

  Their linear span is a complement of $\rad\,F\Pi_A$ in $F\Pi_A$.
\end{corollary}

\begin{proof}
  We have $T\preceq^\dagger S$ if and only if $S=T$, for all $S,T\in
  \Pi_A^<$. Hence the claim follows from Theorem~\ref{multiformel},
  Corollary~\ref{idem-cor} and Theorem~\ref{indecs}.
\end{proof}

\breakitdown

Let $C_A=(C_{P,Q})$ denote the Cartan matrix of $F\Pi_A$, that is,
$C_{P,Q}$ equals the multiplicity of $M_P$ in a composition series of
$\Lambda_Q$, for all $P,Q\in\Pi_A^<$.
Equivalently, 
$$
C_{P,Q}
=
\dim\,e_P\wedge\Lambda_Q
=
\dim\,e_P\wedge F\Pi_A\wedge e_Q\,,
$$
since $M_P$ has dimension one.  We already now that 
\begin{center}
  $C_{P,Q}\neq 0$ implies $Q\preceq P$,
\end{center}
by~\eqref{mult-triangle} and Theorem~\ref{indecs}.  In this case,
Proposition~\ref{lie-help} and Theorem~\ref{multiformel} yield the
following linear basis of the Peirce component 
$e_P\wedge F\Pi_A\wedge e_Q$.

\begin{theorem} \label{cartan}
  Let $P=(P_1,\ldots,P_l),Q=(Q_1,\ldots,Q_k)\in\Pi_A^<$ such that
  $Q\preceq P$, then the elements
  $$
  e_{Q(1)^\circ}\vee\cdots\vee e_{Q(l)^\circ}
  $$
  with $Q(i)\in\Pi_{P_i}^\dagger$ for all $i\in[l]$ and
  $Q(1)\vee \cdots \vee Q(l)\ass Q$, form a linear basis
  of $e_P\wedge F\Pi_A\wedge e_Q$.
  In particular, 
  $$
  C_{P,Q} = (m_1-1)!\cdots (m_l-1)!,
  $$
  where
  $m_j=\#\setofall{i\in[k]}{Q_i\subseteq P_j}$ for all $j\in[l]$.
\end{theorem}

\begin{table}[htbp]
  $$
  \begin{array}{|l||c|c|c|c|c|}
    \hline
    & \rule{0pt}{2.5ex}
      \Lambda_{(123)} 
    & \Lambda_{(12,3)} 
    & \Lambda_{(13,2)} 
    & \Lambda_{(1,23)} 
    & \Lambda_{(1,2,3)} \\[2pt]\hline\hline
    M_{(1,2,3)} &   0   &   0    &    0   &   0    &    1 \\[1mm]\hline
    M_{(1,23)}  &   0   &   0    &    0   &   1    &    1 \\[1mm]\hline
    M_{(13,2)}  &   0   &   0    &    1   &   0    &    1 \\[1mm]\hline
    M_{(12,3)}  &   0   &   1    &    0   &   0    &    1 \\[1mm]\hline
    M_{(123)}   &   1   &   1    &    1   &   1    &    2 \\[1mm]\hline
  \end{array}
  $$
  \caption{Cartan matrix of $F\Pi_3$.}
  \label{tab1}
\end{table}
The Cartan matrix of $F\Pi_3$ is displayed in Table~3. 

\begin{remarks} \label{remarks-cartan} \begin{rm}
  \textbf{(1)} For each $Q\in\Pi_A^<$, the sum of the $Q$-column of
  $C_A$ is $\length{Q}!$, the dimension of $\Lambda_Q$.  The sum of
  the $P$-row of $C_A$ is equal to the dimension of the right ideal
  $e_P\wedge F\Pi_A$ of $F\Pi_A$, for all $P\in\Pi_A^<$. This right
  ideal is indecomposable as an $F\Pi_A$-right module.  If
  $\type{P}=(p_1,\ldots,p_l)$, then there is the explicit formula
  \begin{equation}
    \label{dim-right}
    \dim\,e_P\wedge F\Pi_A=2^l\,|\Pi_{p_1-1}|\,\cdots\,|\Pi_{p_l-1}|\,.
  \end{equation}
  To see this, let $n=|A|$ and consider first the case where $P=(A)$. 
  We have
  $$
  \dim\,e_A\wedge F\Pi_A
  =
  \sum_{Q\in\Pi_A^<} C_{(A),Q}
  =
  \sum_{Q\in\Pi_A^<} (\length{Q}-1)!
  =
  \#\Pi_A^*
  =
  \#\Pi_n^*.
  $$
	To any $R=(R_1,\ldots,R_l)\in\Pi_{\{2,\ldots,n\}}$
	can be associated two set compositions
	$(\{1\},R_1,\ldots,R_l)$ and 
	$(R_1\cup\{1\},R_2,\ldots,R_l)$ in $\Pi_n^*$.
	The identity $\#\Pi_n^*=2\#\Pi_{n-1}$ readily follows, completing
	the proof of~\eqref{dim-right} in this case.
	For example, we have $\dim\,e_{(123)}F\Pi_3=2\#\Pi_2=6$.  For
  arbitrary $P\in\Pi_A^<$, \eqref{dim-right} now follows from the
  factorisation
  $$
  \dim\,e_P\wedge F\Pi_A
  =
  \sum_{Q\in\Pi_A^<} C_{P,Q}
  = 
  \sum_{\substack{Q\in\Pi_A^<\\[2pt] Q\preceq P}} C_{P,Q}
  = 
  \prod_{i=1}^l \sum_{Q\in\Pi_{P_i}^<} C_{P_i,Q}\,.
  $$
   
  \textbf{(2)} The Cartan matrices $C_A$, $A\in\Fin$, have a very
  simple \emph{fractal structure}. Namely, if
  $Q=(Q_1,\ldots,Q_k)\in\Pi^<_A$ is of length $k$, then the non-zero
  part of the $Q$-column of $C_A$ coincides with the (rightmost)
  $(1,2,\ldots,k)$-column of $C_{[k]}$ when $Q_i$ is ``identified''
  with $i$ for all $i\in[k]$.  More formally, let
  $P=(P_1,\ldots,P_l)\in\Pi_A^<$ such that $Q\preceq P$, then
  \begin{equation} \label{fractal}
    C_{P,Q}=C_{I,(1,\ldots,k)}\,,
  \end{equation}
  where $I=(I_1,\ldots,I_l)$ is the set composition in $\Pi_k$ such that
  $Q_{I_j}\in\Pi_{P_j}$ for all $j\in[l]$.  (Indeed, both values are
  equal to $(|I_1|-1)!\cdots(|I_l|-1)!$, by Theorem~\ref{cartan}.)
  
  There is the following interesting structural explanation for the
  identity~\eqref{fractal}.
	Let $K$ denote the annihilator of $\Lambda_Q$ in $F\Pi_A$, then 
  $$
  F\Pi_A=K\oplus F\Pi_Q
  $$
	by~\eqref{left-mod}, where $\Pi_Q$ is the sub-semigroup of $\Pi_A$
	consisting of all set partitions $P\in\Pi_A$ such that $Q\preceq P$.
	There is an isomorphism of algebras $\iota:F\Pi_k\to F\Pi_Q$ mapping
  $$
  (I_1,\ldots,I_l)
  \longmapsto 
  \Big(\bigcup_{i\in I_1} Q_i,\bigcup_{i\in I_2} Q_i,
  \ldots,\bigcup_{i\in I_l} Q_i\Big)
  $$
  for all $(I_1,\ldots,I_l)\in\Pi_k$.
  We have $e_I\iota\equiv e_{I\iota}$ modulo $K$ for all $I\in \Pi_k$,
  since $\min\,Q_1<\cdots<\min\,Q_k$.  
  
  Using the isomorphism $\iota$, $\Lambda_Q$ may be viewed as a module
  for $F\Pi_k$. More precisely, the mapping
  $e_{(i_1,\ldots,i_k)}\mapsto e_{(Q_{i_1},\ldots,Q_{i_k})}$
  defines an isomorphism of $F\Pi_k$-modules from
  $\Lambda_{(1,\ldots,k)}$ onto $\Lambda_Q$.  Thus, if $P\in\Pi_Q$ and
  $I\in\Pi_k$ such that $I\iota=P$, then
  $$
  C_{P,Q}
  =
  \dim\,e_{I\iota}\wedge\Lambda_Q
  =
  \dim\,e_I\iota\wedge\Lambda_Q
  =
  \dim\,e_I\wedge\Lambda_{(1,\ldots,k)}
  =
  C_{I,(1,\ldots,k)},
  $$
  recovering~\eqref{fractal}.     
\end{rm}\end{remarks}

\section{Descending Loewy series} \label{7}
Let $\N_0=\N\cup\{0\}$.  Suppose $\Alg$ is a finite-dimensional
associative algebra over $F$ with identity.  If $M$ is an
$\Alg$-module, then the \emph{descending Loewy series}
$(\rad_\Alg^{(k)}M)_{k\in\N_0}$ of $M$ is defined recursively by
$\rad_\Alg^{(0)}M=M$ and
$
\rad_\Alg^{(k+1)}M
:=
\rad_\Alg(\rad_\Alg^{(k)}M)
=
(\rad\,\Alg)(\rad_\Alg^{(k)}M)
$
for all $k\in\N_0$. The \emph{$k$th Loewy layer} of $M$ is
$\rad_\Alg^{(k)}M/\rad_\Alg^{(k+1)}M$. It is the largest semi-simple
quotient of $\rad_\Alg^{(k)}M$ as an $\Alg$-module. If $M=\Alg$ is the
regular $\Alg$-module, we obtain the descending Loewy series
$(\rad^{(k)}\Alg)_{k\in\N_0}=(\rad_\Alg^{(k)}\Alg)_{k\in\N_0}$ of the algebra $\Alg$.

Throughout this section, $A\in\Fin$. We will describe the descending
Loewy series of $F\Pi_A$. For convenience, we set $\Alg:=F\Pi_A$ and
write $fg$ instead of $f\wedge g$ for $f,g\in\Alg$.
Furthermore, for $P,Q\in\Pi_A$, we write $P\occursin Q$ if $Q\preceq
P$ and $\length{Q}-\length{P}=1$.

\begin{lemma} \label{rad-fund}
  Let $Q=(Q_1,\ldots,Q_k)\in\Pi_A^<$ and $\ol{Q}=(Q_k,\ldots,Q_1)$, 
  then
  $$
  \rad_\Alg\Lambda_Q=\sum_{P\occursin Q} \Lambda_Pe_{\ol{Q}}\;.
  $$
\end{lemma}

\begin{proof}
  If $k=1$, then both sides are equal to zero.  Suppose $k>1$. The
  right hand side is a nilpotent left ideal of $\Alg$, by
  \eqref{mult-triangle}, and contained in $\Lambda_Q$, thus actually
  contained in $\rad_{\Alg}\Lambda_Q$.
  
  Let $Q'=(Q'_1,\ldots,Q'_k)\ass Q$ and $i\in[k-1]$.  Set
  $P:=(Q'_1,\ldots,Q'_{i-1},Q'_i\cup Q'_{i+1},Q'_{i+2},\ldots,Q'_k)$,
  then $P\occursin Q$. In fact, $Q\in\Pi_A^<$ implies that
  $\ol{Q}\preceq^\dagger P$. Hence there exists a sign
  $\varepsilon\in\{+1,-1\}$ such that
  $$
  (*)\quad
  \varepsilon\Big(
  e_{(Q'_1,\ldots,Q'_i,Q'_{i+1},\ldots,Q'_k)}
  -
  e_{(Q'_1,\ldots,Q'_{i+1},Q'_i,\ldots,Q'_k)}
  \Big)
  =
  e_Pe_{\ol{Q}}
  \in
  \sum_{P\occursin Q} \Lambda_Pe_{\ol{Q}},
  $$
  by Theorem~\ref{multiformel}.  The elements of the form $(*)$
  linearly generate $\rad_{\Alg}\Lambda_Q$, by Theorem~\ref{indecs}.
  This completes the proof.
\end{proof}

As a consequence of Lemma~\ref{rad-fund}, we have
\begin{equation} \label{fullrad-fund}
\rad\,\Alg
=
\sum_{Q\in\Pi_A^<}\rad_\Alg\Lambda_Q
=
\spanofall{e_Se_T}{S,T\in\Pi_A,\,S\occursin T}_F\,.
\end{equation}
Using the triangularity property~\eqref{mult-triangle}, it is a
routine matter now to obtain a description of
$\rad^{(k)}_\Alg\Lambda_Q$ (and thus of the $k$th Loewy layer of
$\Lambda_Q$), by induction.

\begin{theorem} \label{loewy-compl}
  Let $Q\in\Pi_A^<$ and $k\in\N_0$, then
  $$
  \rad^{(k)}_{\Alg}\Lambda_Q
  =
  \sum_{\substack{Q\preceq R\in\Pi_A\\[3pt] \length{Q}-\length{R}\ge k}} 
  e_R\Lambda_Q
  =
  \rad^{(k+1)}_{\Alg}\Lambda_Q
  \oplus
  \bigoplus_{\substack{Q\preceq T\in\Pi_A^<\\[3pt] \length{Q}-\length{T}=k}} 
  e_T\Lambda_Q\,.
  $$
\end{theorem}

\begin{proof}
  For $k=0$, this follows from Theorem~\ref{indecs}.
  Let $k>0$ and choose $Q'\ass Q$.
  Suppose $R\in\Pi_A$ such that $Q\preceq R$ and $\length{Q}-\length{R}\ge k$. Then
  $e_Re_{Q'}\in\rad_\Alg\Lambda_Q$, in fact,
  $e_Re_{Q'}\in e_R\rad_\Alg\Lambda_Q$. Hence 
  \begin{eqnarray*}
    e_Re_{Q'}
    & \in &
    \spanofall{(e_Re_P)(e_Pe_{\ol{Q}})}{R\succeq P\occursin Q}_F,
    \mbox{ by Lemma~\ref{rad-fund} and~\eqref{mult-triangle}}\\[3mm]
    & \subseteq &
    (\rad^{(k-1)}\Alg)(\rad_{\Alg}\Lambda_Q),
    \mbox{ by induction}\\[3mm]
    & = &
    \rad^{(k)}_{\Alg}\Lambda_Q\,.
  \end{eqnarray*}
  Conversely,
  \begin{eqnarray*}
    \lefteqn{%
      \rad^{(k)}_\Alg \Lambda_Q
      = 
      (\rad\,\Alg)\,\rad^{(k-1)}_\Alg\Lambda_Q}\hspace*{12ex}\\[3mm]
    & \subseteq &
    \sum_{\substack{S\occursin T\succeq R\succeq Q\\[3pt] \length{Q}-\length{R}\ge k-1}}
     (e_Se_T)(e_R\Lambda_Q)
     \subseteq 
    \sum_{\substack{Q\preceq S\in\Pi_A\\[3pt] \length{Q}-\length{S}\ge k}} e_S\Lambda_Q\,,
  \end{eqnarray*}
  by \eqref{mult-triangle}, \eqref{fullrad-fund} and induction.
  This proves the first equality and implies the inclusion from right
  to left in the second equality.
  
  To prove the remaining part of the second equality, choose
  $R\in\Pi_A$ such that $Q\preceq R$ and $\length{Q}-\length{R}\ge k$, then by
  Corollary~\ref{complement} and \eqref{mult-triangle},
  $$
  e_R\Lambda_Q
  =
  \Big(\sum_{T\in\Pi_A^<}e_T\Big)e_R\Lambda_Q
  =
  \sum_{R\preceq T\in\Pi_A^<} 
  e_Te_R\Lambda_Q
  \subseteq
  e_{T'}\Lambda_Q
  +
  \sum_{\substack{Q\preceq T\in\Pi_A^<\\[2pt] \length{Q}-\length{T}\ge k+1}} 
  e_T\Lambda_Q\,,
  $$
  where $T'$ is the set composition in $\Pi_A^<$ with $T'\ass R$.
  Finally, the sum
  $$
  \rad^{(k+1)}_\Alg\Lambda_Q
  +
  \sum_{\substack{Q\preceq T\in\Pi_A^<\\[2pt] \length{Q}-\length{T}=k}} e_T\Lambda_Q 
  $$
  is direct, by~\eqref{mult-triangle}, Corollary~\ref{complement}
  and the description of $\rad^{(k+1)}_\Alg\Lambda_Q$ given by the
  first equality.
\end{proof}

From the second description of $\rad^{(k)}_\Alg\Lambda_Q$ above, it follows that 
$$
e_T\rad^{(k)}_\Alg\Lambda_Q\not\equiv 0
\mbox{ modulo $\rad^{(k+1)}_\Alg\Lambda_Q$}
$$
if and only if $Q\preceq T$ and $\length{Q}-\length{T}=k$, for all
$Q,T\in\Pi_A^<$.  This gives:

\begin{corollary} \label{mod-struct-loewy}
  Let $Q,T\in\Pi_A^<$, then $M_T$ occurs as a direct summand 
  in the $k$th Loewy layer of $\Lambda_Q$ if and only if 
  $Q\preceq T$ and $\length{Q}-\length{T}=k$.
\end{corollary}

Thus the occurrence of $M_T$ in $\Lambda_Q$ is restricted to a
\emph{single} Loewy layer (if it occurs at all, that is, if $Q\preceq
T$), with multiplicity $C_{T,Q}$.  For example, the Loewy structure of
$\Lambda_{(1,2,3)}$ may be obtained by inspection of the last column
of Table~3. It is displayed in Figure~\ref{fig1}.
\begin{figure}[htbp]
  \centering
$$
\xymatrix@C=-3ex@R=1.5ex{%
                                              && M_{(1,2,3)}   && \\
\Lambda_{(1,2,3)} \quad  =  \quad M_{(12,3)}  && M_{(13,2)}    && M_{(1,23)}\\
                                        & M_{(123)}   &&   M_{(123)}   &  
}
$$
  \caption{Loewy structure of $\Lambda_{(1,2,3)}$}
  \label{fig1}
  \medskip
\end{figure}

As another consequence, we can determine the nil\-index of 
$\rad\,F\Pi_A$, that is, the smallest $k$ such that $\rad^{(k)}F\Pi_A=0$.

\begin{corollary} \label{nilindex}
  The nilindex of $\rad\,F\Pi_A$ is equal to $|A|$.
  More precisely, if $A=\{a_1,\ldots,a_n\}$ is of order $n$ and 
  $a_1<\ldots<a_n$, then 
  \begin{eqnarray*}
    \lefteqn{%
      \rad^{(n-1)}F\Pi_A}\\[3mm]
    & = &
    e_A\wedge F\Pi_A\wedge e_{(a_1,a_2,\ldots,a_n)}
    =
    \spanofall{e_{(x_1,\ldots,x_n)^\circ}}{%
      (x_1,\ldots,x_n)\in\Pi_A,\,x_n=a_1}_F
  \end{eqnarray*}
  has dimension $(n-1)!$, while $ \rad^{(n)}F\Pi_A = 0$.
\end{corollary}

\begin{proof}
  If $P,Q\in\Pi^<_A$, then $Q\preceq P$ implies $0\le \length{Q}-\length{P}\le
  n-1$, with equality on the right if and only if
  $Q=(a_1,a_2,\ldots,a_n)$ and $P=(A)$.  Hence the claims follow from
  Theorems~\ref{loewy-compl} and~\ref{cartan}.
\end{proof}

Let $\Pi^\dagger=\bigcup_{A\in\Fin} \Pi_A^\dagger$.  To conlude this
section, we construct a linear basis of $F\Pi_A$ which is well adapted
to the descending Loewy series of $F\Pi_A$:

\begin{corollary} \label{loewy-basis}
  The elements $e_{Q(1)^\circ}\vee\cdots\vee e_{Q(m)^\circ}$, where
  \vspace*{-2ex}
  \begin{enumerate}
  \item[(a)] $m\in\N_0$;
  \item[(b)] $Q(i)\in\Pi^\dagger$ for all $i\in[m]$;
  \item[(c)] $Q=Q(1)\vee\cdots\vee Q(m)\in\Pi_A$;
  \item[(d)] $\Big(\bigcup Q(1),\ldots,\bigcup Q(l)\Big)\in\Pi_A^<$;
  \end{enumerate}
  form a linear basis of $F\Pi_A$.  If $k\in\N_0$, then those amongst
  these basis elements for which 
  \begin{enumerate}
  \item[(e)] $\length{Q}=m+k$,
  \end{enumerate}
  span a linear complement of $\rad^{(k+1)}F\Pi_A$ in
  $\rad^{(k)}F\Pi_A$.
\end{corollary}

\begin{proof}
  For $k\in\N_0$, we have
  $$
  \rad^{(k)}F\Pi_A
  =
  \bigoplus_{Q\in\Pi_A^<} \rad^{(k)}_{F\Pi_A}\Lambda_Q
  =
  \rad^{(k+1)} F\Pi_A
  \oplus
  \bigoplus_{Q\in\Pi_A^<}
  \bigoplus_{\substack{Q\preceq T\in\Pi_A^<\\[3pt] \length{Q}-\length{T}=k}} 
  e_T\wedge\Lambda_Q
  $$
  by Theorem~\ref{indecs} and the second equality in
  Theorem~\ref{loewy-compl}. Applying Theorem~\ref{cartan} to each of
  the spaces $e_T\wedge\Lambda_Q$ on the right, shows that those
  elements satisfying all five conditions (a)--(e) are linearly
  independent and span a linear complement of $\rad^{(k+1)}F\Pi_A$ in
  $\rad^{(k)}F\Pi_A$, as asserted.
\end{proof}

\begin{remark}
  Corollary~\ref{loewy-basis} has a
  transparent explanation in terms of a Poin\-ca\-r\'e-Birkhoff-Witt
  theorem for $(F\Pi,\vee,\Delta)$, viewed as a free $\Fin$-graded
  bialgebra, with one generator in each degree.

  Indeed, the primitive Lie algebra of $F\Pi$ is generated by the
  elements $e_{Q^\circ}$, $Q\in\Pi^\dagger$, by
  Corollary~\ref{primlie}, Theorem~\ref{multiformel} and
  Proposition~\ref{lie-help}.  (These are the basis elements for $m=1$
  above.)  Corollary~\ref{loewy-basis} says that a basis of $F\Pi$
  is obtained by taking (non-vanishing) increasing concatenation
  products of these basis elements of $\Prim\,F\Pi$, with respect to a
  certain order on the basis.
  Products of basis elements $b_1\in F\Pi_A\cap\Prim\,F\Pi$,
  $b_2\in F\Pi_B\cap\Prim\,F\Pi$ are zero unless $A$ and $B$ are
  disjoint. It is therefore enough to consider this case and to define
  $b_1<b_2$ if $\min\,A<\min\,B$. This imposes condition (d) above. 
  
  From this point of view, the parameter $m$ given in (a) is the 
  number of
  Lie factors in the PBW basis element, condition (c) picks out
  basis elements of total degree $A$, and condition (b) ensures 
  that basis elements of $\Prim\,F\Pi$ occur as factors only.
  
  It seems quite remarkable that there is such a simple connection
  between the PBW structure of the external algebra $(F\Pi,\vee)$ and
  the module structure of the internal algebra $(F\Pi,\wedge)$.  The
  same phenomenon for the Solomon algebra was discovered by Blessenohl
  and Laue~\cite[Theorem~2.1]{blelau02}.  It is not hard to derive
  from the above the following analogue of their result:
  
  \emph{%
    For each $k\in\N_0$, we have
    $\rad^{(k)}(F\Pi,\wedge)=\gamma^{(k)}(F\Pi,\vee)$, where
    $\gamma^{(k)}(F\Pi,\vee)$ is the $k$th component of the descending
    central series of $(F\Pi,\vee)$.}
    
  (For details on the descending central series,
  see~\cite[Section~2]{blelau02}.)
\end{remark}

\section{Ext-quiver} \label{8}
Let $A$ be a finite subset of $\N$, and recall that we write
$P\occursin Q$ if $P\preceq Q$ and $\length{Q}-\length{P}=1$, that is,
if $P$ is obtained by assembling two components of $Q$ (and keeping
the remaining ones, in some order).  
Bearing in mind the considerations at the end of
Section~\ref{2}, we see that $(\Pi_A^<,\occursin)$ is isomorphic to the
Hasse diagram of the support lattice of $(\Pi_A,\wedge_A)$.  As a
consequence, the Ext-quiver of $F\Pi_A$ (which encodes the occurrence
of an irreducible module $M_P$ in the first Loewy layer of a principal
indecomposable module $\Lambda_Q$, see~\cite[4.1.6]{benson98}) has the
following simple description.

\begin{theorem} \label{quiver}
  The Ext-quiver of $F\Pi_A$ is given by the Hasse diagram of the
  support lattice of $(\Pi_A,\wedge_A)$, that is, it is isomorphic 
  to $(\Pi_A^<,\occursin)$.
\end{theorem}

\begin{proof}
  For all $P,Q\in\Pi_A^<$, the irreducible module $M_P$ occurs in 
  the first Loewy layer of $\Lambda_Q$ if and only if $P\occursin Q$, 
  by Corollary~\ref{mod-struct-loewy}, with multiplicity $C_{P,Q}=1$, 
  by Theorem~\ref{cartan}.
\end{proof}

The Ext-quiver of $F\Pi_3$ is displayed in Figure~\ref{q3}.
\begin{figure}[htbp]
  \centering
$$
\xymatrix{%
                      & (\{1\},\{2\},\{3\})  & \\
  (\{1,2\},\{3\})
  \ar@{->}[ur]        & (\{1,3\},\{2\})
                        \ar@{->}[u]            & (\{1\},\{2,3\})
                                                 \ar@{->}[ul]\\
                      & (\{1,2,3\})  
                        \ar@{->}[ul]
                        \ar@{->}[u]
                        \ar@{->}[ur]           &  
}
$$
  \caption{$(\Pi_3^<,\occursin)$}
  \label{q3}
\end{figure}

Theorem~\ref{quiver} suggests that there might be a link between the
Ext-quiver of $FS$ and the support lattice of $S$ for a larger class
of idempotent semigroups~$S$; for example, for those semigroups
associated to arbitrary Coxeter complexes. In fact,
Theorem~\ref{quiver} remains true when dihedral groups (instead of
symmetric groups) are considered.

\section{The module structure of the Solomon algebra} \label{9}
Let $n\in\N$. The Solomon algebra $\DD_n$ is isomorphic to the ring of
$S_n$-invariants, $\BB_n$, of the Solomon-Tits algebra $\Z\Pi_n$, as
explained in Section~\ref{3}. In this section, we use this isomorphism to
relate the module structure of $F\Pi_n$ to the module structure of
$\DD_{n,F}$.  Here we assume that $F$ is \emph{a field of
  characteristic zero}. (The modular case is not yet understood.)  
We will drop the index $F$ in our notations and write, for instance,
$\DD_n$ instead of $\DD_{n,F}$ in what follows.

\subsection{Lie idempotents}
We start with a brief analysis of $\Delta$-primitive idempotents in 
the Solomon algebra.

For each finite subset $A$ of $\N$, let $S_A$ denote the symmetric
group on $A$. If $A$ has order $n$, then $S_A$ is isomorphic to $S_n$,
and everything we have said in Section~\ref{3} remains true when $A$ instead of $[n]$ is considered. In particular, the algebra of 
$S_A$-invariants $\BB_A$ in $F\Pi_A$ is a subalgebra of $F\Pi_A$ 
and isomorphic to the Solomon algebra $\DD_A$ of $S_A$.  The linear map
$$
f\longmapsto \ol{f}=\sum_{\pi\in S_A} f^\pi
$$
is (up to the factor $\frac{1}{n!}$) a projection from $F\Pi_A$ onto
$\BB_A$, since $F$ has characteristic zero.

If $\pi\in S_A$ and $X\subseteq A$, then the natural
action of $\pi$ on the subsets of $A$ gives rise to a linear map
$F\Pi_X\to F\Pi_{X\pi}$, defined by~\eqref{action}, which we also
denote by~$\pi$.  It is easy to see that
\begin{equation} \label{act-coproduct}
 \Delta(f^\pi)=\Delta(f)^{\pi\tensor\pi}
\end{equation}
for all $f\in F\Pi_A$. For, it suffices to consider $f=P\in\Pi_A$,
by linearity. We have
$(P^\pi)|_X=(P|_{X\pi^{-1}})^\pi$ for all $X\subseteq A$, hence
$$
\Delta(P^\pi)
=
\sum_{X\subseteq A} 
\Big(P|_{X\pi^{-1}}\tensor P|_{A\backslash X\pi^{-1}}\Big)^{\pi\tensor\pi}
=
\Delta(P)^{\pi\tensor\pi}.
$$
Combining~\eqref{act-coproduct} with Corollary~\ref{cor-rp}, we obtain
the following result.

\begin{proposition} \label{primitiv-sym}
  Suppose $A\in\Fin$ has order $n$.  If $E$ is a $\Delta$-primitive
  idempotent in $F\Pi_A$, then $\frac{1}{n!}\,\ol{E}$ is a 
  $\Delta$-primitive idempotent in $\BB_A$.
\end{proposition}

The Solomon algebra is intimately linked to the free Lie algebra, as
was discovered by Garsia and Reutenauer~\cite{garsia-reutenauer89}.
The $\Delta$-primitive idempotents in~$\BB_n$ (or, via
Theorem~\ref{th-bidigare}, in $\DD_n$) are easily identified as the
Lie idempotents in~$\BB_n$.  To recall their definition, consider the
free associative algebra $\Alg(X)=\bigoplus_{n\ge 0} \Alg_n(X)$ over an
infinite set $X$.  The $n$th homogeneous component, $\Alg_n(X)$, of
$\Alg(X)$ is an $FS_n$-module, via Polya action:
$$
\pi\,x_1\cdots x_n
:=
x_{1\pi}\cdots x_{n\pi}\,,
$$
for all $\pi\in S_n$, $x_1,\ldots,x_n\in X$. The Lie commutator
$a\circ b=ab-ba$ defines the structure of a Lie algebra on $\Alg(X)$.
Let $\LL(X)=\bigoplus_{n\ge 1} \LL_n(X)$ denote the Lie subalgebra of
$\Alg(X)$ generated by $X$. Then $\LL(X)$ is a free Lie algebra, freely
generated by $X$.  The (right-normed)
Dynkin operator $\omega_n\in FS_n$ can be defined by
$
\omega_n(x_1\cdots x_n)
=
x_1\circ(x_2\circ(x_3\circ\;\cdots\;(x_{n-1}\circ x_n)\cdots))
$
for all $x_1,\ldots,x_n\in X$. The Dynkin-Specht-Wever theorem~\cite{dynkin47,specht48,wever49} says
that $\omega_n^2=n\omega_n$ or, equivalently, that left action
of $\frac{1}{n}\omega_n$ yields a projection from $\Alg_n(X)$ 
onto $\LL_n(X)$. Any such idempotent $e$ in
$FS_n$ is a \emph{Lie idempotent}.  Equivalently, $e\in FS_n$ is a Lie
idempotent if and only if $e\omega_n=\omega_n$ and $\omega_ne=ne$.

If $q=(q_1,\ldots,q_k)$ is a composition of $n$, let 
$\length{q}=k$ and $q^*:=q_1$.  It
is a simple, yet striking observation that
$$
\omega_n
=
\sum_{q\comp n} (-1)^{\length{q}-1} q^* \Xi^q,
$$
hence that $\omega_n\in\DD_n$ (see~\cite[Proposition~1.2]{blelau96}
for the left-normed version of this identity).  
The corresponding element in $\BB_n$ is
$$
\Omega_n
=
\sum_{q\comp n} (-1)^{\length{q}-1} q^* X^q.
$$
Accordingly, we shall refer to elements $E$ in $\BB_n$ such that
$\Omega_n\wedge E=nE$ and $E\wedge\Omega_n=\Omega_n$ as Lie
idempotents in $\BB_n$.
From the definition of the idempotent $e_{[n]}$ given in
Lemma~\ref{primitiv}, it is readily seen that
$$
\Omega_n
=
\frac{1}{(n-1)!}\,\ol{e_{[n]}}.
$$
Hence the idempotent $e_{[n]}$ we used for the study of the
Solomon-Tits algebra is a refinement of the classical Dynkin
operator, and Lemma~\ref{primitiv} is a refinement of the
Dynkin-Specht-Wever theorem, by Proposition~\ref{primitiv-sym}.  
These refined idempotents are defined over
$\Z$ (unlike Lie idempotents). This allowed us to study the
Solomon-Tits algebra over a field of arbitrary characteristic.

\begin{proposition}
  The $\Delta$-primitive idempotents in $\BB_n$ are precisely the Lie
  idempotents in $\BB_n$.
\end{proposition}

This is closely related to~\cite[Theorem~3.1]{parisII}. However,
some care must be taken when linking the coproduct on the direct sum $\bigoplus_n\DD_n$ considered in~\cite{parisII} to the coproduct
$\Delta$ (for more details, see~\cite[Lemma~30]{patras-schocker03}).

\begin{proof}
  Let $e$ be a $\Delta$-primitive idempotent in $\BB_n$.
  From Corollary~\ref{cor-rp} it follows that
  $X^q\wedge e=\sum_{\type{Q}=q}Q\wedge e=0$ whenever $q\comp n$ 
  with $\length{q}>1$, hence that $\Omega_n\wedge e=ne$. 
  Similarly, we get
  $e\wedge\Omega_n=\Omega_n$, since $\frac{1}{n}\Omega_n$ is
  also a $\Delta$-primitive idempotent in $\BB_n$
  by Proposition~\ref{primitiv-sym}. Hence $e$ is a Lie idempotent 
  in $\BB_n$.
  
  Conversely, if $e$ is a Lie idempotent in $\BB_n$, then
  $e=\frac{1}{n}\,\Omega_n\wedge e$. Hence $e$ is $\Delta$-primitive,
  by Propositions~\ref{bi-wedge} and~\ref{primitiv-sym}.
\end{proof}

\subsection{Principal indecomposable modules}  
Let $n\in\N$, and let $\Omega_A:=\frac{1}{|A|!}\ol{e_A}$ denote the
Dynkin operator in $\BB_A$ (up to the factor $|A|$), for all 
$A\in\Fin$. Each $\Omega_A$ is a $\Delta$-primitive
idempotent in $\BB_A$ by Proposition~\ref{primitiv-sym}.  The
corresponding linear basis $\setofall{\Omega_Q}{Q\in\Pi_n}$ of $F\Pi_n$
consists of primitive idempotents such that, additionally,
\begin{equation}
  \label{idem-action}
  \Omega_Q{}^\pi=\Omega_{Q^\pi}
\end{equation}
for all $Q\in\Pi_n$, $\pi\in S_n$.

Suppose $Q\in\Pi_n$ has type $q=(q_1,\ldots,q_k)$, 
then the order of the stabiliser of $Q$ in $S_n$ is 
$s_q=q_1!\cdots q_k!$, and
\begin{equation}
  \label{hli}  
  \ol{\Omega}_Q
  =
  \sum_{\pi\in S_n} \Omega_{Q^\pi}
  =
  s_q\sum_{\type{Q'}=q} \Omega_{Q'}
\end{equation}
depends on $q$ only. We set $f_q:=\ol{\Omega}_Q$. Then the principal
left ideal
$$
\Lambda_q:=F\Pi_n\wedge f_q
$$
of $F\Pi_n$ is $S_n$-invariant and contains $\Omega_P$ whenever
$\type{P}\ass q$, since
\begin{equation}
  \label{sym-multi1}
  \Omega_P\wedge f_q
  =
  s_q
  \sum_{\type{Q'}=q} \Omega_P\wedge \Omega_{Q'}
  =
  s_qc_q \Omega_P
\end{equation}
in this case, by \eqref{prim-idem}. Here $c_q$ is the coefficient
defined in the proof of Corollary~\ref{rad-Dn}. Combined
with~\eqref{hli} and Theorem~\ref{indecs}, this implies that
$$
\Lambda_q
=
\bigoplus_{\substack{Q\in\Pi_n^<\\[2pt] \type{Q}\ass q}} \Lambda_Q
=
\spanofall{\Omega_P}{\type{P}\ass q}_F\,,
$$
where $\Lambda_Q=F\Pi_n\wedge\Omega_Q$ for all $Q\in\Pi_n^<$.
Hence there is the $S_n$-invariant decomposition
$$
F\Pi_n=\bigoplus_{p\p n} \Lambda_p
$$
of $F\Pi_n$ into left ideals by Theorem~\ref{indecs}. These two
formulae combined with \eqref{idem-action}, \eqref{sym-multi1} and
Corollary~\ref{rad-Dn}, allow us to recover the following
result of Garsia and Reutenauer~\cite{garsia-reutenauer89}
and Blessenohl and Laue~\cite{blelau96}.

\begin{theorem}
  Let $n\in\N$, then
  $$
  \BB_n=\bigoplus_{p\p n} \ol{\Lambda}_p
  $$
  is a decomposition of $\BB_n$ into indecomposable left ideals.
  For each $p\p n$, $\ol{\Lambda}_p$ has $F$-basis
  $\setofall{f_q}{q\ass p}$ consisting (up to non-zero scalar factors)
  of primitive idempotents, and $\rad_{\BB_n}\ol{\Lambda}_p$ has
  $F$-basis $\{f_p-f_q\,|\,q\ass p,\,q\neq p\}$.
  
  Furthermore, the modules 
  $M_p=\ol{\Lambda}_p/\rad_{\BB_n}\ol{\Lambda}_p$, $p\p n$, 
  have dimension one and form a complete set of irreducible
  $\BB_n$-modules.
\end{theorem}

For later use, note that if $Q=(Q_1,\ldots,Q_k)\in\Pi_n$ has 
type $q$ and $S^Q$ is a right transversal of the stabiliser of $Q$ 
in $S_n$, then
\begin{equation} \label{e-symm}
	\ol{e_{\raisebox{-2pt}{$\scriptstyle Q$}}}
	=
	\sum_{\sigma\in S^Q}
	\Big(
	\ol{e_{\raisebox{-2pt}{$\scriptstyle Q_1$}}}
	\vee\cdots\vee
	\ol{e_{\raisebox{-2pt}{$\scriptstyle Q_k$}}}
	\Big)^\sigma
	=
	s_q\sum_{\sigma\in S^Q}\Omega_Q{}^\sigma
	=
	\ol{\Omega}_Q
	=
	f_q\,.
\end{equation}

\subsection{%
	The Solomon-Tits algebra as a module\\ for the Solomon algebra}
Let $n\in\N$. It is natural to consider $F\Pi_n$ as a module 
for the subalgebra~$\BB_n$. 
We will study the action of $\BB_n$ on the basis $\setofall{e_Q}{Q\in\Pi_n}$ and the 
indecomposable modules $\Lambda_Q=F\Pi_n\wedge e_Q$ of $F\Pi_n$.

\begin{proposition} \label{Bn-iso}
  Let $Q,R\in\Pi_n^<$, then $\Lambda_Q$ and $\Lambda_R$ are 
  isomorphic as $\BB_n$-modules if and only if $\type{Q}\ass\type{R}$.
  Furthermore, $M_Q$ is an irreducible $\BB_n$-module isomorphic to
  $M_p$, where $p\p n$ such that $\type{Q}\ass p$.
\end{proposition}

\begin{proof}
  It is more convenient here to consider the Dynkin operator 
  $\Omega_A$ instead of $e_A$ for all $A\in\Fin$. We may
  do so by Remark~\ref{idem-choice}.

  If $\type{Q}\ass\type{R}$, then there exist $\tilde{Q}\ass Q$ 
  and $\pi\in S_n$ such that $\tilde{Q}^\pi=R$. The mapping
  $\Lambda_Q\to\Lambda_R$ which
  sends $f\mapsto f^\pi$ is then an isomorphism of $\BB_n$-modules,
  by~\eqref{idem-action}.
  Conversely, if $\Lambda_Q$ and $\Lambda_R$ are isomorphic as
  $\BB_n$-modules, 
  then $f_q\wedge \Lambda_Q\neq 0$ implies $f_q\wedge \Lambda_R\neq
  0$, hence $\tilde{R}\preceq Q^\pi$ for some $\tilde{R}\ass R$, 
  $\pi\in S_n$,
  by~\eqref{hli} and~\eqref{mult-triangle}. Interchanging the roles of
  $Q$ and $R$, there exist $\tilde{Q}\ass Q$ and $\sigma\in S_n$ 
  such that
  $\tilde{Q}\preceq R^\sigma$. This implies $\type{Q}\ass \type{R}$,
  completing the proof of the first claim.
  
  The second claim also follows from~\eqref{hli} 
  and~\eqref{mult-triangle}, since
  \[
  f_p\wedge \Omega_Q
  =
  s_p\sum_{\type{P}=p} \Omega_P\wedge \Omega_Q
  =
  s_p\sum_{\substack{\type{P}=p\\[2pt] P\ass Q}} \Omega_P
  \notin
  \rad_{F\Pi_n}\Lambda_Q\,.\qedhere
  \]
\end{proof}

Let $q\comp n$ and $Q=(Q_1,\ldots,Q_k)\in\Pi_n$ be of type $q$, then
$\pi\in S_n$ is a \emph{block permutation} of $Q$ if $\pi$ performs a
blockwise permutation of the components (blocks) of $Q$, that is,
more precisely, if $Q^\pi\ass Q$ and $\pi|_{Q_i}$ is non-decreasing for
all $i\in[k]$.  For instance, $\pi=2\,1\,6\,5\,4\,3\in S_6$ is a block
permutation of $Q=(13,5,4,26)$, while $\sigma=6\,1\,2\,5\,4\,3\in S_6$
is not a block permutation of $Q$ (although
$Q^\pi=(26,4,5,13)=Q^\sigma$).  Let $G_Q$ denote the subgroup of $S_n$
consisting of all block permutations of $Q$ and set
$\alpha_Q=1/|G_Q|\sum_{\pi\in G_Q} \pi$.

\begin{lemma} \label{splatter}
	Let $Q,\tilde{Q},Q'\in\Pi_n$ such that $Q\ass\tilde{Q}\ass Q'$.
	Then, for all $\pi\in G_Q$, we have 
	$
	e_{\tilde{Q}}{}^\pi
	=
	e_{\tilde{Q}^\pi}
	$.
	Furthermore, $e_{\tilde{Q}}{}^{\alpha_Q}=e_{Q'}{}^{\alpha_Q}$
	if and only if $\type{\tilde{Q}}=\type{Q'}$.
\end{lemma}

\begin{proof}
	If $A,B\in\Fin$ and $\tau:A\to B$ is an order-preserving bijection,
	then $e_A{}^\tau=e_B$.
	Hence
	$$
	e_{\tilde{Q}}{}^\pi
	=
	e_{\tilde{Q}_1}{}^\pi\vee\cdots\vee e_{\tilde{Q}_k}{}^\pi
	=
	e_{\tilde{Q}_1\pi}\vee\cdots\vee e_{\tilde{Q}_k\pi}
	=
	e_{\tilde{Q}^\pi}
	$$
	for all $\pi\in G_Q$. Furthermore, $\type{\tilde{Q}}=\type{Q'}$
	if and only if there exists a block permutation $\pi\in G_Q$
	such that $\tilde{Q}^\pi=Q'$. By the above, this is equivalent
	to $e_{\tilde{Q}}{}^{\alpha_Q}=e_{Q'}{}^{\alpha_Q}$.
\end{proof}

As a consequence, the indecomposable module
$\Lambda_Q=\spanofall{e_{\tilde{Q}}}{\tilde{Q}\ass Q}_F$
is invariant under the action of $G_Q$ and therefore a 
$(\BB_n,FG_Q)$-bimodule.

\begin{proposition} \label{abspalter}
  Suppose $Q\in\Pi_n$ is of type $q$ and $p\p n$ such that 
  $q\ass p$, then
  $
  \Lambda_Q
  =
  \Lambda_Q{}^{\alpha_Q}\oplus\Lambda_Q{}^{\mathrm{id}_n-\alpha_Q}
  $
  is a decomposition into $\BB_n$-submodules, and
  $\Lambda_Q{}^{\alpha_Q}$ is isomorphic to the indecomposable 
  $\BB_n$-module $\ol{\Lambda}_p$.
\end{proposition}

Here $\mathrm{id}_n$ denotes the identity in $S_n$.

\begin{proof}
  The first part is obvious. For the proof of the second part, 
  choose a right transversal $G^Q$ of $G_Q$ in $S_n$.
  Then the map
	$\varphi:\Lambda_Q{}^{\alpha_Q}\to\ol{\Lambda}_p$ 
	which sends 
	$g\mapsto \sum_{\pi\in G^Q}g^\pi$ is an isomorphism 
	of $\BB_n$-modules, since
	$e_{\tilde{Q}}{}^{\alpha_Q}\varphi=1/|G_Q|\,f_{\tilde{q}}$
	for all $\tilde{Q}\ass Q$ of type $\tilde{q}$, 
	by~\eqref{e-symm}.
\end{proof}

There is a surprising link here to a fundamental result of Garsia and
Reute\-nau\-er~\cite[Theorem~4.5]{garsia-reutenauer89} on the Solomon
algebra.  This result states that an element $\gamma\in FS_n$ lies in
the Solomon algebra $\DD_n$ if and only if, via Polya action,
$$
\gamma\, l_1\cdots l_k 
\in 
\spanofall{l_{1\pi}\cdots l_{k\pi}}{\pi\in S_k}_F\,,
$$
for all $k$-tuples of homogeneous Lie elements $l=(l_1,\ldots,l_k)$
such that $l_1\cdots l_k\in \Alg_n(X)$.  As a consequence, there is an
action of $\DD_n$ on the linear span $U_l$ of the elements
$l_{1\pi}\cdots l_{k\pi}$ ($\pi\in S_k$). The modules $U_l$
are characteristic for the Solomon algebra in the sense that $\DD_n$
is the common stabiliser in $FS_n$ of these linear spaces with
respect to Polya action. In fact, it suffices 
to consider \emph{generic} tuples $l=(l_1,\ldots,l_k)$ of Lie
monomials only for which the elements $l_{1\pi}\cdots l_{k\pi}$
($\pi\in S_k$) are linearly independent. The \emph{type} of $l$ is the
composition $q=(q_1,\ldots,q_k)$ of $n$ such that $l_i\in \LL_{q_i}(X)$
for all $i\in[k]$.

The $\DD_n$-modules $U_l$ may, of course, be viewed as modules 
for $\BB_n$.

\begin{theorem} \label{char-iso}
  Let $n\in\N$ and $q\comp n$. Suppose that $l=(l_1,\ldots,l_k)$ is a
  generic $k$-tuple of homogeneous Lie monomials, of type $q$.  Then
  $U_l\cong \Lambda_Q$ as a $\BB_n$-module, for any $Q\in\Pi_n$ of
  type $q$.
\end{theorem}

\begin{proof}
  If $I=\{i_1,\ldots,i_m\}\subseteq [k]$ such that $i_1<\cdots<i_m$, put
  $
  q_I := (q_{i_1},\ldots,q_{i_k})
  $
  and 
  $
  l_I := l_{i_1}\cdots l_{i_k}
  $.
  Let $r=(r_1,\ldots,r_m)\comp n$.  Then, on the one hand,
  $$
  \Xi^rl_1\cdots l_k 
  = 
  \sum_{%
    \substack{(I_1,\ldots, I_m)\in\Pi_k\\[2pt]  q_{I_j}\comp r_j}}
  l_{I_1}\cdots l_{I_m}\,,
  $$
  by~\cite[Theorem~2.1]{garsia-reutenauer89}.  If
  $Q=(Q_1,\ldots,Q_k)\in\Pi_n$ has type $q$, then, on the other hand,
  there is $1$-$1$ correspondence between the set compositions
  $I=(I_1,\ldots,I_m)\in\Pi_k$ such that $q_{I_j}\comp r_j$ and the
  set compositions $R\in\Pi_n$ of type $r$ such that $Q\preceq R$,
  namely the map
  $$
  I
  \longmapsto
  R_I
  :=
  \Big(\bigcup_{i\in I_1} Q_i,\bigcup_{i\in I_2} Q_i,
  \ldots,\bigcup_{i\in I_m} Q_i\Big)
  $$
  already considered in Remark~\ref{remarks-cartan}~(2).
  It follows from \eqref{left-mod} that
  \begin{eqnarray*}
    X^r\wedge e_Q
    &= &
    \sum_{%
      \substack{I=(I_1,\ldots, I_m)\in\Pi_k\\[2pt]  q_{I_j}\comp r_j}}
    e_{R_I\wedge Q}
     = 
    \sum_{%
      \substack{(I_1,\ldots, I_m)\in\Pi_k\\[2pt]  q_{I_j}\comp r_j}}
    e_{Q_{I_1}}\vee\cdots\vee e_{Q_{I_m}}\,.
  \end{eqnarray*}
  Thus the map $e_{(Q_{1\pi},\ldots,Q_{k\pi})}\mapsto
  l_{1\pi}\cdots l_{k\pi}$ ($\pi\in S_k$) defines an isomorphism of
  $\BB_n$-modules from $\Lambda_Q$ onto $U_l$.
\end{proof}

Combining Theorem~\ref{char-iso} with Proposition~\ref{abspalter}, we
see that the projective indecomposable $\DD_n$-module $\ol{\Lambda}_p$
is isomorphic to a direct summand of $U_l$ whenever $l$ is of type
$p$. This result is due to Mielck~\cite{mielck96}.
Furthermore, the following description of the decomposition numbers of
the characteristic $\DD_n$-modules $U_l$ is now immediate from
Theorem~\ref{char-iso} and Proposition~\ref{Bn-iso}.

\begin{corollary}
  Let $n\in\N$ and $q\comp n$, and suppose $Q\in\Pi_n^<$ has 
  type $q$.
  Then, for any generic tuple $l$ of
  homogeneous Lie monomials of type $q$, the multiplicity of
  $M_p$, $p\p n$, in a composition series of $U_l$ is equal to
  the sum of all Cartan invariants $C_{P,Q}$ of $F\Pi_n$, 
  taken over set compositions
  $P\in\Pi_n^<$ with $\type{P}\ass p$ .
\end{corollary}

A better understanding of the $\DD_n$-modules $U_l$
would be very interesting. This includes as a (crucial)
special case the study of $FS_n$ as $\DD_n$-module (when $l$ is of
type $(1,1,\ldots,1)$).

\subsection{Cartan invariants}  
The most efficient way to recover the well-known description of the
Cartan invariants of $\DD_n$ (see~\cite{garsia-reutenauer89},\cite[Corollary~2.1]{blelau96},\cite[Section~3.6]{parisII}) seems to build on
Proposition~\ref{abspalter}, as follows.

Let $c_{r,q}$ denote the multiplicity of $M_r$ in a composition 
series of $\ol{\Lambda}_q$, for all $q,r\p n$. 
Suppose $Q\in\Pi_n^<$ has type $\ass q$.
Then, for $\Lambda_Q=F\Pi_n\wedge e_Q$, we obtain
$$
c_{r,q}
=
\dim f_r\wedge\ol{\Lambda}_q
=
\dim f_r\wedge\Lambda_Q{}^{\alpha_Q} 
=
\dim (f_r\wedge\Lambda_Q)^{\alpha_Q}.
$$
We will need the following auxiliary result.

\begin{lemma}
  Let $r\p n$, $Q\in\Pi_n^<$, and set
  $$
	X
	:=
	\bigoplus_{%
  		\substack{Q\preceq T\in\Pi_n^<\\[2pt] \type{T}\ass r}}
  		e_T\wedge\Lambda_Q\,.    
	$$
  Then 
  $f_r\wedge\Lambda_Q=f_r\wedge X$ 
  and 
  $\dim\,f_r\wedge X=\dim\,X$.
\end{lemma}

\begin{proof}
  Set $k:=\length{Q}-\length{r}$. If $k<0$, then
  $f_r\wedge\Lambda_Q=0$ as is readily seen 
  from~\eqref{triangular}, \eqref{left-mod} and~\eqref{hli},
  	while $X=0$ is clear.
  
	Suppose $k\ge 0$, then	any occurrence of $M_R$ in $\Lambda_Q$,
	for $R$ of type $\ass r$, is restricted to the $k$th Loewy layer of 
	$\Lambda_Q$ (if it occurs at all), by 
	Corollary~\ref{mod-struct-loewy}. Hence
  $$
  f_r\wedge\Lambda_Q/\rad^{(k)}_{F\Pi_n}\Lambda_Q
  =
  0
  =
  f_r\wedge\rad^{(k+1)}_{F\Pi_n}\Lambda_Q\,,
  $$
  by Proposition~\ref{Bn-iso} and~\eqref{hli}. 
  Now Theorem~\ref{loewy-compl} implies
  $$
	f_r\wedge\Lambda_Q
	\subseteq
	f_r\wedge(\rad^{(k)}\Lambda_Q)
	=
	f_r\wedge
	\bigoplus_{%
		\substack{Q\preceq T\in\Pi_n^<\\[2pt] \length{Q}-\length{T}=k}}
		e_T\wedge\Lambda_Q\,.    
  $$
  Furthermore, if $T\in\Pi_n^<$ such that $\length{Q}-\length{T}=k$ 
  and $\type{T}\not\ass r$, then $f_r\wedge e_T=0$.
  It follows that $f_r\wedge\Lambda_Q\subseteq f_r\wedge X$, 
  the other inclusion is clear.
  
  The dimensions of $f_r\wedge\Lambda_Q=f_r\wedge X$ and
  of $X$ both describe the
  multiplicity of $M_r$ in a $\BB_n$-composition series of
  $\Lambda_Q$, by Proposition~\ref{Bn-iso}.
\end{proof}

As a consequence of the preceding result, we have
$$
c_{rq}
=
\dim\,(f_r\wedge\Lambda_Q)^{\alpha_Q}
=
\dim\,(f_r\wedge X)^{\alpha_Q}
=
\dim\,X^{\alpha_Q}.
$$
Put $l:=\length{r}$. Then $X^{\alpha_Q}$ has set of linear generators
$(e_{Q(1)^\circ}\vee\cdots\vee e_{Q(l)^\circ})^{\alpha_Q}$
where $Q(i)\in\Pi$ for all $i\in[l]$ such that 
$Q(1)\vee\cdots\vee Q(l)\ass Q$ and
\begin{equation} \label{star-1}
	\Big(\bigcup Q(1),\ldots,\bigcup Q(l)\Big)\in \Pi_n^<
\end{equation} 
has type $\ass r$. Using a triangularity argument, we can drop 
condition~\eqref{star-1} and obtain the set of generators
$$
\sum_{\pi\in S_l} 
(e_{Q(1\pi)^\circ}\vee\cdots\vee e_{Q(l\pi)^\circ})^{\alpha_Q}.
$$
From Lemma~\ref{splatter} we know that the map 
which sends $e_{\tilde{Q}}{}^{\alpha_Q}\mapsto \type{\tilde{Q}}$ is a linear isomorphism
from $\Lambda_Q{}^{\alpha_Q}$ onto the $F$-linear span 
(in the free associative algebra $\Alg(\N)$ over $\N$) of 
$\setofall{\tilde{q}}{\tilde{q}\ass q}$.
It maps $X^{\alpha_Q}$ onto the linear span $T_{r,q}$ of the elements
$$
\sum_{\pi\in S_l} q(1\pi)^\circ.\,\ldots\,.q(l\pi)^\circ
$$
where $q(i)\comp r_i$ for all $i\in[l]$ such that 
$q(1).\,\ldots\,.q(l)\ass q$.  Here $r.s$ is the (concatenation) product of the compositions $r$ and $s$ in $\Alg(\N)$ and 
$q\mapsto q^\circ$ is the associated right-normed Dynkin mapping.
Thus $c_{r,q}=\dim\,T_{r,q}$.

Let $B$ be a set of compositions and $\beta:B\to \LL(\N)$ be an injective
mapping such that $B\beta$ is a linear basis of $\LL(\N)$.
Consider a total order on $B$ such that $q\le r$ whenever 
the sum of the components of $q$ is less than the sum of the components of $r$.
The corresponding (symmetrised) 
Poincar\'e-Birkhoff-Witt basis of $\Alg(\N)$ consists of 
the elements
\begin{equation} \label{star-2}
	\sum_{\pi\in S_m} (q(1\pi)\beta).\,\ldots\,.(q(m\pi)\beta)
\end{equation}
where $q(i)\in B$ for all $i\in[m]$ and $q(i)\ge q(i+1)$ for all
$i\in[m-1]$.
Those basis elements~\eqref{star-2} for which $m=l$, $q(i)\comp r_i$
for all $i\in[l]$ 
and $q(1).\,\ldots\,.q(l)\ass q$, form a linear basis of
$T_{r,q}$.
Hence $c_{r,q}$ is equal to the number of these basis elements,
in accordance with the results in~\cite{blelau02,garsia-reutenauer89}.
An explicit formula for $c_{r,q}$ follows easily
(see~\cite[Corollary~2.1]{blelau02}).

\subsection{Descending Loewy series}
As a general observation on idempotent semigroups, we know from
Corollary~\ref{cor-brown} that $\rad\,\BB_n=\BB_n\cap\rad\,F\Pi_n$.
Comparison of the results of Section~\ref{7} with those on
the Loewy structure of $\DD_n$ derived in \cite{blelau96} yields the
following surprising generalisation.

\begin{theorem} \label{loewy-fix}
  $\rad^{(k)}\BB_n=\BB_n\cap\rad^{(k)}F\Pi_n$, for all $k\in\N_0$.
\end{theorem}

\begin{proof}
	For $k\in\N_0$, $\rad^{(k)}F\Pi_n$ has set of linear generators
	$$
	e_{\raisebox{-1pt}{$\scriptstyle Q(1)^\circ$}}
	\vee\cdots\vee 
	e_{\raisebox{-1pt}{$\scriptstyle Q(m)^\circ$}}
	$$
	by Corollary~\ref{loewy-basis}, where $Q(i)\in\Pi$ for all $i\in[m]$
	such that $Q(1)\vee\cdots\vee Q(m)\in\Pi_n$ and 
	$\sum_{i=1}^m\length{Q(i)}\ge m+k$.
	Hence by~\eqref{e-symm}
	$\ol{\rad^{(k)}F\Pi_n}=\BB_n\cap\rad^{(k)}F\Pi_n$
	has set of linear generators 
	$f_{q(1)^\circ.\,\ldots\,.q(m)^\circ}$,
	where $q(i)$ is a composition for all $i\in[m]$
	such that $q(1).\,\ldots\,.q(m)\comp n$ and 
	$\sum_{i=1}^m\length{q(i)}\ge m+k$.
	(The symbol $f$ is understood to be linear with respect 
	to subscripts.)
	These elements linearly generate $\rad^{(k)}\BB_n$,
	by~\cite[Theorem~2.5]{blelau96}.
\end{proof}

Again it would be interesting to know whether 
Theorem~\ref{loewy-fix} holds for other Coxeter groups as well.

\newcommand{\appears}[1]{To appear in \emph{#1}}

\end{document}